%% file: agt-1-6.tex
\input gtmacros

\input amsnames
\input amstex
\let\cal\Cal          
\catcode`\@=12   
%
\input agtout

\volumenumber{1}
\volumeyear{2001}
\papernumber{6}
\published{24 February  2001}
\pagenumbers{115}{141}
\received{29 October 2000}
\revised{16 February 2001}
\accepted{16 February 2001}
\font\bboldss=cmmib10
\def\boldxi{\hbox{\bboldss\char'030}}
\def\bolddelta{\hbox{\bboldss\char'016}}
\let\\\par

\def\heading#1\endheading{{\def\S##1{\relax}\def\\{\relax\ignorespaces}
    \section{#1}}}
\def\head#1\endhead{\heading#1\endheading}

\def\subhead#1\endsubhead{\sh{#1}}
\def\subsubhead#1\endsubsubhead{\sh{#1}}
\def\specialhead#1\endspecialhead{\sh{#1}}
\def\demo#1{\rk{#1}\ignorespaces}
\def\enddemo{\ppar}

\def\qed{\ifmmode\quad\sq\else\hbox{}\hfill$\sq$\par\goodbreak\rm\fi}  
\def\proclaim#1{\rk{#1}\sl\ignorespaces}
\def\endproclaim{\rm\ppar}
\def\cite#1{[#1]}
\newcount\itemnumber

\let\itemold\item
\def\item{\itemold{{\rm(\number\itemnumber)}}%
\global\advance\itemnumber by 1\ignorespaces}
\def\S{section~\ignorespaces}  
\def\SS{sections~\ignorespaces} 
\def\date#1\enddate{\relax}
\def\thanks#1\endthanks{\relax}   
\def\dedicatory#1\enddedicatory{\relax}  
\def\rom#1{{\rm #1}}  
\let\footnote\plainfootnote

\def\Refs{\ppar{\large\bf References}\ppar\bgroup\leftskip=45pt
\frenchspacing\parskip=3pt plus2pt\small}       
\def\endRefs{\egroup}
\def\widestnumber#1#2{\relax}
\def\endrefitem{}
\def\refdef#1#2#3{\def#1{\leavevmode\unskip\endrefitem#2\def\endrefitem{#3}}}
\def\ref{\par}
\def\endref{\endrefitem\par\def\endrefitem{}}
\refdef\key{\noindent\llap\bgroup[}{]\ \ \egroup}
\refdef\no{\noindent\llap\bgroup[}{]\ \ \egroup}
\refdef\by{\bf}{\rm, }
\refdef\manyby{\bf}{\rm, }
\refdef\paper{\it}{\rm, }
\refdef\book{\it}{\rm, }
\refdef\jour{}{ }
\refdef\vol{}{ }
\refdef\yr{(}{) }
\refdef\ed{(}{, Editor) }
\refdef\publ{}{ }
\refdef\inbook{from: ``}{'', }
\refdef\pages{}{ }
\refdef\page{}{ }
\refdef\paperinfo{}{ }
\refdef\bookinfo{}{ }
\refdef\publaddr{}{ }
\refdef\moreref{}{ }
\refdef\eds{(}{, Editors)}
\refdef\bysame{\hbox to 3 em{\hrulefill}\thinspace,}{ }
\refdef\toappear{(to appear)}{ }
\refdef\issue{no.\ }{ }

\title{Generalized Orbifold Euler Characteristic\\of Symmetric 
Products\\\vglue 1pt\\and Equivariant Morava K-Theory}
\covertitle{Generalized Orbifold Euler Characteristic\\of Symmetric 
Products\\\noexpand\vglue 1pt\\and Equivariant Morava K-Theory}
\asciititle{Generalized Orbifold Euler Characteristic\\of Symmetric 
Products\\and Equivariant Morava K-Theory}

\shorttitle{Generalized Orbifold Euler Characteristic}

\author{Hirotaka Tamanoi}

\address{
Department of Mathematics, University of California Santa Cruz,\\ 
Santa Cruz, CA 95064, USA}

\email{tamanoi@math.ucsc.edu}

\keywords{
Equivariant Morava K-theory, generating functions, $G$-sets, M\"obius
functions, orbifold Euler characteristics, q-series, second quantized
manifolds, symmetric products, twisted iterated free loop space,
twisted mapping space, wreath products, Riemann zeta function
}
\asciikeywords{
Equivariant Morava K-theory, generating functions, G-sets, Moebius
functions, orbifold Euler characteristics, q-series, second quantized
manifolds, symmetric products, twisted iterated free loop space,
twisted mapping space, wreath products, Riemann zeta function
}

\primaryclass{55N20, 55N91}
\secondaryclass{57S17, 57D15, 20E22, 37F20, 05A15}

\abstract
We introduce the notion of generalized orbifold Euler characteristic
associated to an arbitrary group, and study its properties. We then
calculate generating functions of higher order ($p$-primary) orbifold
Euler characteristic of symmetric products of a $G$-manifold $M$. As a
corollary, we obtain a formula for the number of conjugacy classes of
$d$-tuples of mutually commuting elements (of order powers of $p$) in
the wreath product $G\wr\frak S_n$ in terms of corresponding numbers
of $G$. As a topological application, we present generating functions
of Euler characteristic of equivariant Morava K-theories of symmetric
products of a $G$-manifold $M$.
\endabstract

\asciiabstract{
We introduce the notion of generalized orbifold Euler characteristic
associated to an arbitrary group, and study its properties. We then
calculate generating functions of higher order (p-primary) orbifold
Euler characteristic of symmetric products of a G-manifold M. As a
corollary, we obtain a formula for the number of conjugacy classes of
d-tuples of mutually commuting elements (of order powers of p) in
the wreath product G wreath S_n in terms of corresponding numbers
of G. As a topological application, we present generating functions
of Euler characteristic of equivariant Morava K-theories of symmetric
products of a G-manifold M.}

\makeshorttitle

\catcode`\@=\active

\head
\S 1 Introduction and summary of results \quad 
\endhead

Let $G$ be a finite group and let $M$ be a smooth $G$-manifold. We
study generalized orbifold Euler characteristics of $(M;G)$. These are
integer-valued invariants associated to any group $K$. (See (1-3)
below.)  The simplest of such invariants (corresponding to the trivial
group $K=\{e\}$) is the usual Euler characteristic $\chi(M/G)$ of the
orbit space. It is well known that $\chi(M/G)$ can be calculated as
the average over $g\in G$ of Euler characteristic of corresponding
fixed point submanifolds:
$$
\chi(M/G)=\frac1{|G|}\sum_{g\in G}\chi(M^{\langle g\rangle}),
\tag1-1
$$
where $\langle g\rangle\le G$ is the subgroup generated by $g\in G$. 
See for example, \cite{Sh, p.127}. 

In 1980s, string physicists proposed a notion of {\it orbifold Euler
characteristic} of $(M;G)$ defined by
$$
\chi_{\text{orb}}(M;G)=\frac1{|G|}\sum_{gh=hg}\chi(M^{\langle
g,h\rangle}),
\tag1-2
$$
where the summation is over all commuting pairs of elements in $G$
\cite{DHVW}. The orbifold Euler characteristic is always an integer,
since (1-1) implies 
$$
\chi_{\text{orb}}(M;G)=\sum_{[g]}\chi(M^{\langle g\rangle}/C_G(g))\in
\Bbb Z,
$$
where the summation is over all the conjugacy classes of $G$, and
$C_G(g)$ is the centralizer of $g$ in $G$. This formula is of the
form $\chi(M/G)+(\text{correction terms})$. 

\subhead
Generalized orbifold Euler characteristic
\endsubhead
Let $K$ be any group. The generalized orbifold Euler characteristic of
$(M;G)$ associated to $K$ is an integer 
$$
\chi_{\sssize K}(M;G)\overset{\text{def}}\to=\!\!\!\!\!\!\!\!\!\!\!\!\!
\sum_{[\phi]\in\text{Hom}(K,G)/G}\!\!\!\!\!\!\!\!\!\!\!\!\!
\chi\bigl(M^{\langle \phi\rangle}/C_G(\phi)\bigr)
=\dfrac1{|G|}\!\!\!\!\!\!\!\!\!\!\!\!\!\!\!\!\!\!\!\!
\sum_{\ \ \ \ \ \ \ \ \ \ \phi\in\text{Hom}(K\times\Bbb Z,G)}
\!\!\!\!\!\!\!\!\!\!\!\!\!\!\!\!\!\!\!\!\!\!\chi(M^{\langle \phi\rangle}).
\tag1-3
$$
The first summation is over $G$-conjugacy classes of homomorphisms,
and the second equality is a consequence of (1-1). Thus, either
expression can be taken as the definition of $\chi_K(M;G)$. Here,
$C_G(\phi)$ is the centralizer in $G$ of the image
$\langle\phi\rangle$ of $\phi$. Note that when $K$ is the trivial
group $\{e\}$ or $\Bbb Z$, our $\chi_K(M;G)$ specializes to (1-1) or
(1-2). In \S2, we describe its various properties including
multiplicativity and the following formula for a product $K\times L$
of two groups:
$$
\chi_{\sssize K\times L}(M;G)=\!\!\!\!\!\!\!\!\!\!\!\!\!
\sum_{[\phi]\in\text{Hom}(K,G)/G}\!\!\!\!\!\!\!\!\!\!\!\!\!
\chi_{\sssize L}\bigr(M^{\langle\phi\rangle};C_G(\phi)\bigr).
\tag1-4
$$
This formula, which is easy to prove, is crucial for inductive steps
in the proofs of our main results, Theorems A and B below.

In this paper, we are mostly concerned with the cases $K=\Bbb Z^d$ and
$K=\Bbb Z^d_p$, where $\Bbb Z_p$ denotes the ring of $p$-adic integers. We
use the following notations:
$$
\chi_{\Bbb Z^d}(M;G)=\chi^{(d)}(M;G),\qquad
\chi_{\Bbb Z^d_p}(M;G)=\chi^{(d)}_p(M;G).
\tag1-5
$$
We call these $d$-th order ($p$-primary) orbifold Euler
characteristics. 

Our definition (1-3) is partly motivated by consideration of a mapping
space $\text{Map}(\Sigma, M/G)$, where a manifold $\Sigma$ has
fundamental group $K$. When $\Sigma$ is the genus $g$ orientable
surface $\Sigma_g$ with $\pi_1(\Sigma_g)=\Gamma_g$, we call the
corresponding quantity $\chi_{\Gamma_g}(M;G)$ genus $g$ orbifold Euler
characteristic of $(M;G)$.

\subhead
Higher order orbifold Euler characteristics for symmetric products
\endsubhead
It turns out that $\chi^{(d)}(M;G)$ admits a geometric interpretation
in terms of the mapping space $\text{Map}(T^d,M/G)$, where $T^d$ is
the $d$-dimensional torus. See \S 2 for more details. And as such, it
is very well behaved. We demonstrate this point by calculating
$\chi^{(d)}$ for symmetric products.

When $M$ is a $G$-manifold, the $n$-fold Cartesian product $M^n$
admits an action of a wreath product $G\wr\frak S_n$. The orbit space
$M^n/(G\wr\frak S_n)=SP^n(M/G)$ is the $n$-th symmetric product
of $M/G$.

\proclaim{Theorem A} For any $d\ge0$ and for any $G$-manifold $M$, 
$$
\sum_{n\ge0}q^n\chi^{(d)}(M^n;G\wr\frak S_n)
=\Bigl[\prod_{r\ge1}(1-q^r)^{j_r(\Bbb
Z^d)}\Bigr]^{(-1)\chi^{(d)}(M;G)},
\tag1-6
$$
where $j_r(\Bbb Z^d)=\!\!\!\!\!\!\!\!\sum\limits_{r_1r_2\cdots r_d=r}
\!\!\!\!\!\!\!\!r_2r_3^2\cdots r_d^{d-1}$ is the number of index
$r$ subgroups in $\Bbb Z^d$. 
\endproclaim

For the case of $d=1$ ($\chi^{(1)}=\chi_{\text{orb}}$), the above
formula was calculated by Wang \cite{W}. By letting $M$ to be a point,
and using the notation $G_n=G\wr\frak S_n$, we obtain

\proclaim{Corollary 1-1} For any $d\ge0$, we have 
$$
\sum_{n\ge0}q^n\bigl|\text{\rm Hom}\bigl(\Bbb Z^d, G_n\bigr)/G_n\bigr|
=\Bigl[\prod_{r\ge1}(1-q^r)^{j_r(\Bbb Z^d)}\Bigr]^{(-1)
|\text{Hom}(\Bbb Z^d,G)/G|}.
\tag1-7
$$
\endproclaim

Special cases of our results are known. Macdonald \cite{M1} calculated
Euler characteristic of symmetric products of any topological space
$X$. His formula reads
$$
\sum_{n\ge0}q^n\chi\bigl(SP^n(X)\bigr)
=\frac1{(1-q)^{\chi(X)}}.
\tag1-8
$$
Hirzebruch-H\"ofer \cite{HH} calculated the orbifold Euler
characteristic (1-2), which is our $\chi^{(1)}(M;G)$, of symmetric
products. Their formula is
$$
\sum_{n\ge0}q^n\chi_{\text{orb}}(M^n;\frak S_n)
=\Bigl[\prod_{r\ge1}(1-q^r)\Bigr]^{(-1)\chi(M)}.
\tag1-9
$$

After completing this research, the author became aware of the paper
\cite{BF} in which the formulae (1-6) and (1-7) with trivial $G$ were
calculated. (That the exponent of their formula can be identified with
$j_r(\Bbb Z^d)$ was pointed out by Allan Edmonds in Math\. Review.) We
remark that the formula (1-7) with trivial $G$ is straightforward once
we observe that $\bigl|\text{Hom}(\Bbb Z^d,\frak S_n)/\frak S_n\bigr|$
is the number of isomorphism classes of $\Bbb Z^d$-sets of order $n$,
and $j_r(\Bbb Z^d)$ is the number of isomorphism classes of transitive
(irreducible) $\Bbb Z^d$-sets of order $r$. The second fact is because
the isotropy subgroup of transitive $\Bbb Z^d$-sets of order $r$ is a
sublattice of index $r$ in $\Bbb Z^d$. See also an exercise and its
solution in \cite{St, p.76, p.113}. On the other hand, when $G$ is
nontrivial, a geometric interpretation of elements in $\text{Hom}(\Bbb
Z^d,G_n)/G_n$ is more involved. Our method of proving (1-6) is a
systematic use of the formula (1-4) for generalized orbifold Euler
characteristic and the knowledge of centralizers of elements of the
wreath product $G_n$ described in detail in \S 3. Our method can also
be applied to more general context including $p$-primary orbifold
Euler characteristic $\chi^{(d)}_p(M;G)$.

The integer $j_r(\Bbb Z^d)$ has very interesting number theoretic
properties. It is easy to prove that the Dirichlet series whose
coefficients are $j_r(\Bbb Z^d)$ can be expressed as a product of
Riemann zeta functions $\zeta(s)=\sum_{n\ge1}1/n^s$ with $s\in\Bbb C$: 
$$
\sum_{n\ge1}\dfrac{j_n(\Bbb Z^d)}{n^s}
=\zeta(s)\zeta(s-1)\cdots\zeta(s-d+1).
\tag1-10
$$
For the history of this result, see \cite{So}. 

\subhead
Euler characteristic of equivariant Morava K-theory of symmetric
products
\endsubhead
Let $K(d)^*(X)$ be the $d$-th Morava K-theory of $X$ for
$d\ge0$. Since $K(d)^*$ is a graded field, we can count the dimension
of $K(d)^*(X)$ over $K(d)^*$, if it is finite. We are interested in
computing Euler characteristic of equivariant $d$-th Morava $K$-theory
of a $G$-manifold $M$:
$$
\chi\bigl(K(d)^*_G(M)\bigr)=\dim K(d)^{\text{even}}
(EG\times_G M)-\dim K(d)^{\text{odd}}
(EG\times_G M).
\tag1-11
$$
In \cite{HKR}, they calculate this number in terms of M\"obius
functions \cite{HKR, Theorem B}. It is a simple observation to
identify (1-11) as the $d$-th order $p$-primary orbifold Euler
characteristic $\chi^{(d)}_p(M;G)$ (see a paragraph before Proposition
5-1). Our second main result is as follows.

\proclaim{Theorem B} Let $d\ge0$ and let $M$ be a $G$-manifold. The
Euler characteristic of equivariant Morava $K$-theory is equal to the
$d$-th order $p$-primary orbifold Euler characteristic of $(M;G)$\rom{:}
$$
\chi\bigl(K(d)^*_G(M)\bigr)=\chi^{(d)}_p(M;G).
\tag1-12
$$
The generating function of Euler characteristic of equivariant $d$-th
Morava $K$-theory of symmetric products is given by
$$
\sum_{n\ge0}q^n\chi\bigl(K(d)^*_{G_n}(M^n)\bigr)
=\Bigl[\prod_{r\ge0}(1-q^{p^r})^{j_{p^r}(\Bbb Z^d_p)}\Bigr]^{(-1)
\chi(K(d)^*_G(M))}.
\tag1-13
$$
Here, $G_n=G\wr\frak S_n$, and $j_{\ell}(\Bbb Z^d_p)$ is the number of
index $\ell$ subgroups in $\Bbb Z^d_p$ given by
$$
j_{p^r}(\Bbb Z^d_p)=\!\!\!\!\!\!\!\!\!
\sum_{r_1r_2\cdots r_d=p^r}\!\!\!\!\!\!\!\!r_2r_3^2\cdots
r_d^{d-1},\ \text{and}\ 
j_{\ell}(\Bbb Z^d_p)=0 \text{\ if $\ell$ is not a power of $p$}. 
\tag1-14
$$
\endproclaim

Let $M$ be a point. The resulting formula is both topological and
combinatorial:

\proclaim{Corollary 1-2} For any $d\ge0$, we have 
$$
\aligned
\sum_{n\ge0}q^n\chi\bigl(K(d)^*(BG_n)\bigr)
&=\Bigl[\prod_{r\ge0}(1-q^{p^r})^{j_{p^r}(\Bbb Z^d_p)}
\Bigr]^{(-1)\chi(K(d)^*BG)} \\
\sum_{n\ge0}q^n\bigl|\text{\rm Hom}
(\Bbb Z^d_p,G_n)/G_n\bigr|
&=\Bigl[\prod_{r\ge0}(1-q^{p^r})^{j_{p^r}(\Bbb Z^d_p)}
\Bigr]^{(-1)|\text{Hom}(\Bbb Z^d_p,G)/G|}.
\endaligned
\tag1-15
$$
\endproclaim

When $G$ is a trivial group and hence $G_n=\frak S_n$, the above
formula is straightforward by an argument in terms of $\Bbb
Z_p^d$-sets. 

Again the integers $j_{p^r}(\Bbb Z^d_p)$ have number theoretic
properties and it is well known that the corresponding Dirichlet
series can be expressed as a product of $p$-local factors of Riemann
zeta function. Namely, letting
$\zeta_p(s)=(1-p^{-s})^{(-1)}$ denote the $p$-local factor in the
Euler decomposition of $\zeta(s)$, we have
$$
\sum_{r\ge0}\frac{j_{p^r}(\Bbb Z^d_p)}{p^{rs}}
=\zeta_p(s)\zeta_p(s-1)\cdots\zeta_p(s-d+1).
\tag1-16
$$
In particular, we have the following Euler decomposition of Dirichlet
series:
$$
\sum_{n\ge1}\frac{j_n(\Bbb Z^d)}{n^s}
=\!\!\!\!\!\prod_{p:\text{prime}}\!\!
\biggl(\sum_{r\ge0}\frac{j_{p^r}(\Bbb Z^d_p)}{p^{rs}}\biggr).
\tag1-17
$$

In \cite{H}, Hopkins shows that when $G$ is a trivial group, the
exponents in (1-15) satisfy (1-16) by a method completely different
from ours: by integrating a certain function over $\text{GL}_n(\Bbb
Q_p)$. And he identifies these exponents as Gaussian binomial
coefficients \cite{M2, p.26}.

There is a physical reason why symmetric products of manifolds give
rise to very interesting generating functions. In physics, the process
of quantization of the state space of particles or strings moving on a
manifold produces a Hilbert space of quantum states. The second
quantization then corresponds to taking the total symmetric products
of this Hilbert space, describing quantum states of many particles or
strings. Reversing the order of these two procedures, we can first
apply the second quantization to the manifold, by taking the total
symmetric products of the manifold. This object quantizes well, and in
\cite{DMVV} (see also \cite{D}) they calculate complex elliptic genera 
of second quantized K\"ahler manifolds and it is shown that they are
genus $2$ Siegel modular forms whose weight depend on holomorphic
Euler characteristic of K\"ahler manifolds.

The organization of this paper is as follows. In section 2, we define
generalized orbifold Euler characteristics associated to arbitrary
groups, and describe their properties. We also express abelian
orbifold Euler characteristics in terms of M\"obius functions. We also
describe a geometry behind our definition of generalized orbifold
Euler characteristic. In section 3, we describe some properties of
wreath products including the structure of centralizers. Material here
is not new. However, detailed description on this topic seems to be
rather hard to find in literature, so we worked out details and we
decided to include it. This section is purely group theoretic and is
independent from the rest of the paper. In sections 4 and 5, we
compute higher order ($p$-primary) orbifold Euler characteristics of
symmetric products in the form of generating functions. A relation to
Euler characteristic of equivariant Morava $K$-theory is discussed in
section 5.

{\it Acknowledgement.} The author thanks M\. J\. Hopkins for useful
discussions during the author's visit at MIT, and for making his
preprint \cite{H} available. Our Theorem B was worked out after seeing
his paper. The author also thanks N\. J\. Kuhn who, during the initial
circulation of this paper, informed him that in \cite{K} he had
considered a fixed point functor $F_K$ on $G$-CW complexes which is
formally related to our generalized orbifold Euler characteristic.

\head
\S 2 Generalized orbifold Euler characteristics
\endhead

A generalization of physicists' orbifold Euler characteristic (1-2)
was given in the introduction in (1-3). Here, group $K$ can be an
arbitrary group. Properties enjoyed by $\chi^{(d)}(M;G)$ and
$\chi^{(d)}_p(M;G)$ given in (1-5) become transparent in this
generality. If the group $K$ is abelian, our generalized orbifold
Euler characteristic is better behaved and it admits an expression in
terms of M\"obius functions defined on the family of abelian subgroups
of $G$.

Later in this section, we will explain a geometric meaning of
generalized orbifold Euler characteristics in terms of twisted mapping
spaces.

\subhead
Generalized orbifold Euler characteristics
\endsubhead
We prove basic properties of generalized orbifold Euler
characteristic. Recall from \S 1 that the generalized orbifold Euler
characteristic associated to a group $K$ is given by
$$
\chi_{\sssize K}(M;G)
=\frac1{|G|}\!\!\!\!\!\!\!\!\!\!\!\!\!\!\!\!\!\!\!\!\!\!\!\!\!
\sum_{\ \ \ \ \ \ \ \ \ \ \ \ \phi\in\text{Hom}(K\times\Bbb Z,G)}
\!\!\!\!\!\!\!\!\!\!\!\!\!\!\!\!\!\!\!\!\!\!\!\!\!\chi(M^{\langle \phi\rangle})
=\!\!\!\!\!\!\!\!\!\!\!\!\!\!\!\!\!\!\!\!\!\!\!\!\!\!
\sum_{\ \ \ \ \ \ \ \ \ \ \ \ [\phi]\in\text{Hom}(K,G)/G}
\!\!\!\!\!\!\!\!\!\!\!\!\!\!\!\!\!\!\!\!\!\!\!\!\!\!
\chi\bigl(M^{\langle \phi\rangle}/C_G(\phi)\bigr)
\tag2-1
$$
See the last subsection of \S 2 for a geometric motivation of this
definition. Letting $M$ be a point, we obtain a useful formula
$$
\chi_{\sssize K}(pt;G)=\frac{\bigl|\text{Hom}(K\times\Bbb Z,G)\bigr|}{|G|}
=\bigl|\text{Hom}(K,G)/G\bigr|.
\tag2-2
$$
Now we give a proof of (1-4). 

\proclaim{Proposition 2-1} The orbifold Euler characteristic
$\chi_{\sssize K}$ is multiplicative.\break Namely, for any $G_i$-manifolds
$M_i$ for $i=1,2$, we have 
$$
\chi_{\sssize K}(M_1\times M_2;G_1\times G_2)=
\chi_{\sssize K}(M_1;G_1)\cdot\chi_{\sssize K}(M_2;G_2).
\tag2-3
$$
Furthermore, for any two groups $K$ and $L$, we have 
$$
\chi_{\sssize K\times L}(M;G)=\!\!\!\!\!\!\!\!\!\!\!\!\!
\sum_{[\phi]\in\text{Hom}(K,G)/G}\!\!\!\!\!\!\!\!\!\!\!\!\!
\chi_{\sssize L}\bigr(M^{\langle\phi\rangle};C_G(\phi)\bigr).
\tag2-4
$$
\endproclaim
\demo{Proof} The first formula is straightforward from the definition
of $\chi_K$ given in (2-1). For the second formula, first we observe
that $\chi_{\sssize L}\bigr(M^{\langle\phi\rangle};C_G(\phi)\bigr)$
depends only on the conjugacy class of $\phi$, so the formula is well
defined. Now
$$
\align
(\text{R.H.S})&=\!\!\!\!\!\!\!\!\!\!\!\!\!\!\!\!\!\!
\sum_{\phi\in\text{Hom}(K,G)\ \ \ \ \ \ \ \ }
\!\!\!\!\!\!\!\!\!\!\!\!\!\!\!\!\!
\frac1{\#[\phi]}\frac1{|C_G(\phi)|}
\!\!\!\!\!\!\!\!\!\!\!\!\!\!\!\!\!\!\!\!\!\!\!\!\!\!\!\!
\sum_{\ \ \ \ \ \ \ \ \ \ \ \ \ \ \psi\in\text{Hom}(L\times\Bbb Z, C(\phi))}
\!\!\!\!\!\!\!\!\!\!\!\!\!\!\!\!\!\!\!\!\!\!\!\!\!\!\!\!\!\!
\chi\bigl(M^{\langle \phi,\psi\rangle}\bigr)\\
&=\frac1{|G|}\!\!\!\!\!\!\!\!\!\!\!\!\!\!\!\!\!\!\!\!\!\!\!\!\!\!\!\!
\sum_{\ \ \ \ \ \ \ \ \ \ \ \ \ \ (\phi,\psi)\in\text{Hom}(K\times L\times\Bbb Z,G)}
\!\!\!\!\!\!\!\!\!\!\!\!\!\!\!\!\!\!\!\!\!\!\!\!\!\!\!\!\!\!\!\!\!\!
\chi\bigl(M^{\langle\phi,\psi\rangle}\bigr)
=\chi_{\sssize K\times L}(M;G).
\endalign
$$
Here, $\#[\phi]$ is the number of elements in the conjugacy class of
$\phi$. Since $|C_G(\phi)|$ is the isotropy subgroup of the
conjugation action of $G$ on the homomorphism set at $\phi$, we have
$\#[\phi]\cdot|C_G(\phi)|=|G|$. This completes the proof.
\qed
\enddemo

Next, we rewrite our orbifold Euler characteristic in terms of
M\"obius functions $\mu_H(\,X\,)$ defined for any subgroup $H$ and any
$G$-CW complex $X$. These are defined by downward induction on $P\le
G$ by the formula
$$
\sum_{P\le H\le G}\!\!\!\mu_H(X)=\chi(X^P).
\tag2-5
$$
It is known that any additive functions on $G$-CW complexes can be
expressed as a linear combination of $\mu_H(\ \cdot\ )$'s with $\Bbb
Z[1/|G|]$-coefficients \cite{HKR, Proposition 4.6}. Our generalized
orbifold Euler characteristic has the following expression in terms of
M\"obius functions. 

\proclaim{Lemma 2-2} For any group $K$ and $G$-CW complex $M$, we have
$$
\chi_{\sssize K}(M;G)=\!\!\sum_{H\le G}
\frac{|\text{\rm Hom}(K\times\Bbb Z,H)|}{|G|}
\mu_{\sssize H}(M)
=\!\!\sum_{H\le G}\frac{|H|}{|G|}\cdot\bigl|\text{\rm Hom}(K,H)/H\bigr|\cdot
\mu_{H}(M).
$$
\endproclaim
\demo{Proof}
In the definition of $\chi_{\sssize K}(M;G)$ in (2-1), we replace
$M^{\langle \phi\rangle}$ by (2-5). We have 
$$
\multline
\chi_{\sssize K}(M;G)
=\frac1{|G|}\sum_{\phi}
\chi(M^{\langle \phi\rangle})
=\frac1{|G|}\sum_{\phi}
\sum_{\langle \phi\rangle\le H}\mu_{\sssize H}(M)\\
=\frac1{|G|}\sum_{H\le G}\mu_{\sssize H}(M)
\!\!\!\!\!\!\!\sum_{\phi:K\times\Bbb Z @>>> H}\!\!\!\!\!\!\!1
=\frac1{|G|}\sum_{H\le G}
\bigl|\text{Hom}(K\times\Bbb Z,H)\bigr|\cdot\mu_{\sssize H}(M).
\endmultline
$$
Here in the second and third summation, $\phi$ runs over all
homomorphisms in the set $\text{Hom}(K\times\Bbb Z, G)$. 

The second equality of the statement is due to (2-2).  This proves the
Lemma.
\qed
\enddemo

\subhead
Abelian orbifold Euler characteristics and abelian M\"obius functions
\endsubhead
We recall some facts on complex oriented additive functions 
\cite{HKR, \SS4.1, 4.2}. Let $\theta:\{\text{$G$-CW complexes}\} @>>>
\Bbb Z$ be an integer-valued $G$-homotopy invariant function on $G$-CW
complexes. Then the function $\theta$ is called {\it additive} if it
satisfies
$$
\theta(X\cup Y)+\theta(X\cap Y)=\theta(X)+\theta(Y), \qquad
\theta(\emptyset)=0
$$
for any $G$-CW complexes $X,Y$. For a $G$-equivariant complex
vector bundle $\xi$ on $X$, let $F(\xi)$ be the
associated bundle of complete flags in $\xi$. An additive function
$\theta$ is called {\it complex oriented} if $\theta$ satisfies 
$\theta\bigl(F(\xi)\bigr)=n!\theta(X)$ for any $G$-equivariant complex
$n$-dimensional bundle $\xi$ on $X$. It is known that any complex
oriented additive function on $G$-CW complexes is completely
determined by its value on the family of finite $G$-sets $\{G/A\}$,
where $A$ runs over all abelian subgroups of $G$. In fact, the
following formula holds \cite{HKR, Proposition 4.10}:
$$
\theta(\ \cdot\ )=\frac1{|G|}\!\!\!\!\!
\sum \Sb A\le G \\ A:\text{abelian} \endSb\!\!\!\!\!
|A|\cdot\theta(G/A)\cdot\mu_A^{\Bbb C}(\ \cdot\ ),
\tag2-6
$$
where $\mu_A(\ \cdot\ )$ is a complex oriented additive function
defined by downward induction on an abelian subgroup $A$ by 
$$
\sum\Sb A\le B\le G \\ B:\text{abelian}\endSb
\!\!\!\!\!\mu_B^{\Bbb C}(X)=\chi(X^A)
\tag2-7
$$
for any $G$-CW complex $X$. 

For our generalized orbifold Euler characteristic, when the group $K$
is an abelian group $E$, then $\chi_{\sssize E}(\ \cdot\ ;G)$ is a
complex oriented additive function, since the image of any
homomorphism $\phi:E\times\Bbb Z @>>> G$ is abelian. As such,
$\chi_{\sssize E} (\ \cdot\ ;G)$ satisfies a formula of the form
(2-6). We will explicitly derive this formula in Proposition 2-3. 

On the other hand, we can also show that $\chi_{\sssize E}(\ \cdot\ ;G)$ can
be written as a linear combination of complex oriented additive
functions $\{\chi_{\sssize E}(\ \cdot\ ;A)\}_A$ with $\Bbb
Z[1/|G|]$-coefficients, where $A$ runs over all abelian subgroups of
$G$. For this description, we need a function $\mu_{\!\Cal
A\!}:\{\text{abelian subgroups of $G$}\} @>>> \Bbb Z$ defined by
downward induction on an abelian subgroup $A$ by
$$
\sum \Sb A\le B\le G \\ B:\text{abelian}\endSb \!\!\!\!\!
\mu_{\!\Cal A\!}(B)=1.
\tag2-8
$$
Note that when $G$ is abelian, this relation implies that $\mu_{\!\Cal
A\!}(G)=1$ and $\mu_{\!\Cal A\!}(A)=0$ for any proper subgroup $A$ of
$G$. Thus, (2-8) is of interest only when $G$ is non-abelian. We call
$\mu_A^{\Bbb C}(\ \cdot\ )$ and $\mu_{\Cal A}(\ \cdot\ )$ abelian
M\"obius functions. We rewrite the generalized abelian orbifold
Euler characteristic as follows. In (2-9) below, the first identity
can be proved easily using (2-6), but here we give a different and
amusing proof: we calculate a triple summation in three different
orders. 

\proclaim{Proposition 2-3} Let $E$ be an abelian group. Then the
corresponding orbifold Euler characteristic $\chi_{\sssize E}(M;G)$
satisfies 
$$
\chi_{\sssize E}(M;G)
=\!\sum_{B\le G}\!\frac{|B|}{|G|}\cdot\bigl|\text{\rm Hom}(E,B)\bigr|
\cdot\mu_B^{\Bbb C}(M)
=\!\sum_{A\le G}\!\frac{|A|}{|G|}\cdot\mu_{\Cal A}(A)\cdot
\chi_{\sssize E}(M;A).
\tag2-9
$$
Here in the above summations, $A$ and $B$ run over all abelian
subgroups of $G$. 
\endproclaim
\demo{Proof} We consider the following summation in three variables
$\phi, A,B$: 
$$
(*)=\sum_{\phi, A, B}\mu_{\Cal A}(A)\cdot\mu_B^{\Bbb C}(M),
$$
where $\phi:E\times\Bbb Z @>>> G$ and $A, B$ are abelian subgroups
satisfying $\langle \phi\rangle\le A$ and $\langle \phi\rangle\le
B$. We compute this summation in three different ways: 
$$
(1)\ \sum_A\sum_{\phi}\sum_B,\qquad
(2)\ \sum_{\phi}\sum_A\sum_B,\qquad
(3)\ \sum_B\sum_{\phi}\sum_A.
$$
For the case (1), the summation becomes
$$
\align
(*)&=\sum_{A}\mu_{\Cal A}(A)\!\!\!\!\!\!\!\!\!\!\!\!
\sum_{\phi:E\times\Bbb Z @>>> A\ \ \ \ \ }\!\!\!\!\!\!\!\!\!\!
\Bigl\{\sum_{\langle\phi\rangle\le B}\mu_B^{\Bbb C}(M)\Bigr\}
=\sum_A\mu_{\Cal A}(A)\!\!\!\!\!\!\!\!\sum_{\phi:E\times\Bbb Z @>>> A}
\!\!\!\!\!\!\!\!\chi\bigl(M^{\langle \phi\rangle}\bigr)\\
&=\sum_A\mu_{\Cal A}(A)\cdot|A|\cdot\chi_{\sssize E}(M;A).
\endalign
$$
Here, the summation over $A$ is over all abelian subgroups of $G$,
and (2-7) was used for the second equality. This is allowed since
$\langle\phi\rangle$ is abelian for any homomorphism $\phi:E\times\Bbb
Z @>>> A$. The third equality is the definition of 
$\chi_{\sssize E}(\ \cdot\ ;A)$ in (2-1). 

For the case (2), we have 
$$
(*)=\!\!\!\!\!\!\!\!\!\!\!\sum_{\phi:E\times\Bbb Z @>>> G\ \ \ }
\!\!\!\!\!\!\!\!
\Bigl\{\sum_{\langle\phi\rangle\le A}\mu_{\!\Cal A\!}(A)\Bigr\}
\Bigl\{\sum_{\langle\phi\rangle\le B}\mu_B^{\Bbb C}(M)\Bigr\}
=\!\!\!\!\!\!\!\!\sum_{\phi:E\times\Bbb Z @>>> G}
\!\!\!\!\!\!\!\!1\cdot
\chi\bigl(M^{\langle\phi\rangle}\bigr)
=|G|\cdot\chi_{\sssize E}(M;G).
$$
Note that if $E$ is not abelian, then the image $\langle\phi\rangle$
of $\phi:E\times\Bbb Z @>>> G$ can be non-abelian and the second
equality above may not be valid. This is where we need to assume that
$E$ is abelian. 

For the third summation 
$$
\align
(*)&=\sum_{B}\mu_B^{\Bbb C}(M)\Bigl[\!\!\!\!\!\!\!\!\!\!\!
\sum_{\phi:E\times\Bbb Z @>>> B\ \ \ \ \ }\!\!\!\!\!\!\!\!\!\!\!
\Bigl\{\sum_{\langle\phi\rangle\le A}\mu_{\!\Cal A\!}(A)\Bigr\}\Bigr]
=\sum_B\mu_B^{\Bbb C}(M)\cdot\bigl|\text{Hom}(E\times\Bbb Z,B)\bigr|\\
&=\sum_B|B|\cdot\bigl|\text{Hom}(E,B)\bigr|\cdot\mu_B^{\Bbb C}(M).
\endalign
$$
Here the summation over $B$ is over all abelian subgroups of
$G$. Since $B$ is abelian, $\text{Hom}(E\times \Bbb Z,B)$ is a product
of $\text{Hom}(E,B)$ and $B$. This completes the proof.
\qed
\enddemo

By letting $E$ be the trivial group, we get an interesting formula for
$\chi(M/G)$. 

\proclaim{Corollary 2-4} For any $G$-manifold $M$, we have 
$$
\chi(M/G)=\frac1{|G|}\!\!\!\!\!\!\!\!\!\sum_{\ \ \ \ A:\text{abelian}}
\!\!\!\!\!\!\!\!\!|A|\cdot
\mu_{\!\Cal A\!}(A)\cdot\chi(M/A)
=\frac1{|G|}\!\!\!\!\!\!\!\!\!\sum_{\ \ \ \ B:\text{abelian}}
\!\!\!\!\!\!\!\!\!|B|\cdot\mu_B^{\Bbb C}(M).
$$
\endproclaim
It is interesting to compare this formula with (1-1). Note that the
above formula does not imply that $\mu_{\!\Cal A\!}(A)\cdot\chi(M/A)$
is equal to $\mu_A^{\Bbb C}(M)$. The second equality holds only after
summation over all abelian groups. A similar remark applies to (2-9). 

\subhead
Higher order orbifold Euler characteristic
\endsubhead
We specialize our previous results on generalized orbifold Euler
characteristic to higher order orbifold Euler characteristic 
$\chi^{(d)}(M;G)=\chi_{\Bbb Z^d}(M;G)$. First, (2-1) specializes to 
$$
\chi^{(d)}(M;G)
=\dfrac1{|G|}\!\!\!\!\!\!\!\!\!\!\!\!\!\!\!\!\!\!\!\!\!\!
\sum_{\ \ \ \ \ \ \ \ \ \ \ \phi\in\text{Hom}(\Bbb Z^{d+1},G)}
\!\!\!\!\!\!\!\!\!\!\!\!\!\!\!\!\!\!\!\!\!\!\!\chi(M^{\langle \phi\rangle})
=\!\!\!\!\!\!\!\!\!\!\!\!\!\!\!\!\!\!\!\!\!\!\!\!\!\!
\sum_{\ \ \ \ \ \ \ \ \ \ \ [\phi]\in\text{Hom}(\Bbb Z^d,G)/G}
\!\!\!\!\!\!\!\!\!\!\!\!\!\!\!\!\!\!\!\!\!\!\!\!\!\!
\chi\bigl(M^{\langle \phi\rangle}/C_G(\phi)\bigr).
\tag2-10
$$
Notice that the second equality above is also a consequence of (2-2)
with $K=\Bbb Z^d$ and $L=\{e\}$. If we apply (2-2) with $K=\Bbb Z$ and
$L=\Bbb Z^{d-1}$, then we obtain the following inductive formula. 

\proclaim{Proposition 2-5} For any $d\ge1$, and for any $G$-manifold
$M$, we have 
$$
\chi^{(d)}(M;G)=\sum_{[g]}\chi^{(d-1)}
\bigl(M^{\langle g\rangle};C_G(g)\bigr),
\tag2-11
$$
where the summation is over all conjugacy classes $[g]\in\text{\rm
Hom}(\Bbb Z,G)/G$ of $G$. 
\endproclaim
This is the formula which allows us to prove Theorem A inductively on
$d\ge0$. Lastly, formula (2-9) specializes in our case to 
$$
|G|\cdot\chi^{(d)}(M;G)=\!
\sum_{A\le G}\!|A|\cdot\mu_{\Cal A}(A)\cdot\chi^{(d)}(M;A)
=\!\sum_{B\le G}\!|B|^{d+1}\cdot\mu_B^{\Bbb C}(M).
\tag2-12
$$
Here the summations is over all abelian subgroups of $G$. 

\subhead
Higher order $p$-primary orbifold Euler characteristic
\endsubhead
Recall that the basic formula of this orbifold Euler characteristic
$\chi^{(d)}_p(M;G)=\chi_{\Bbb Z^d_p}(M;G)$ is given by letting 
$K=\Bbb Z^d_p$ in (2-1): 
$$
\chi^{(d)}_p(M;G)
=\frac1{|G|}\!\!\!\!\!\!\!\!\!\!\!\!\!\!\!\!\!\!\!\!\!\!\!\!
\sum_{\ \ \ \ \ \ \ \ \ \ \ \phi\in\text{Hom}(\Bbb Z^d_p\times\Bbb Z,G)}
\!\!\!\!\!\!\!\!\!\!\!\!\!\!\!\!\!\!\!\!\!\!\!\!\chi(M^{\langle \phi\rangle})
=\!\!\!\!\!\!\!\!\!\!\!\!\!\!\!\!\!\!\!\!\!\!\!\!
\sum_{\ \ \ \ \ \ \ \ \ \ [\phi]\in\text{Hom}(\Bbb Z^d_p,G)/G}
\!\!\!\!\!\!\!\!\!\!\!\!\!\!\!\!\!\!\!\!\!\!\!\!
\chi\bigl(M^{\langle \phi\rangle}/C_G(\phi)\bigr).
\tag2-13
$$
Now let $K=\Bbb Z_p$ and $L=\Bbb Z_p^{d-1}$ in the formula (2-4). We
obtain the following inductive formula corresponding to (2-11) for the
$p$-local case. 

\proclaim{Proposition 2-6} For any $d\ge1$ and for any $G$-manifold
$M$, we have 
$$
\chi^{(d)}_p(M;G)=\!\!\!\!\!\!\!\!\!\!\!\!\!\!\!\!\!\!\!\!\!\!\!
\sum_{\ \ \ \ \ \ \ [\phi]\in\text{Hom}(\Bbb Z_p,G)/G}
\!\!\!\!\!\!\!\!\!\!\!\!\!\!\!\!\!\!\!\!\!\!\!
\chi^{(d-1)}_p\bigl(M^{\langle\phi\rangle};C_G(\phi)\bigr).
\tag2-14
$$
Here $[\phi]$ runs over all $G$-conjugacy classes of elements of order
powers of $p$. 
\endproclaim
The following formula, which is a specialization of (2-9) in our
setting, will be used later in \S5 to compare $\chi^{(d)}_p(M;G)$ with
Euler characteristic of equivariant Morava K-theory.

\proclaim{Proposition 2-7} For any $d\ge0$ and for any $G$-manifold
$M$, 
$$
\chi^{(d)}_p(M;G)
=\sum_{A\le G}\frac{|A|}{|G|}\cdot\mu_{\!\Cal A\!}(A)
\cdot\chi^{(d)}_p(M;A)
=\sum_{B\le G}\frac{|B|}{|G|}\cdot |B_{(p)}|^d\cdot
\mu_B^{\Bbb C}(M),
\tag2-15
$$
where the summation is over all abelian subgroups of $G$. 
\endproclaim

\subhead
Generalized orbifold Euler characteristic and twisted mapping space
\endsubhead
We discuss a geometric origin of orbifold Euler
characteristics. Physicists' orbifold Euler characteristic (1-2)
originates in string theory. Higher order ($p$-primary) orbifold Euler
characteristics $\chi^{(d)}(M;G)$ and $\chi^{(d)}_p(M;G)$ have similar
geometric interpretations in terms of twisted mapping spaces. There is
a very strong analogy between this geometric situation and methods
used in orbifold conformal field theory. We can predict results in
orbifold conformal field theory, for example a description of twisted
sectors for the action of wreath products, simply by examining this
geometric situation of twisted mapping spaces.

To describe the geometry, first we consider the free loop space
$L(M/G)=\text{Map}(S^1,M/G)$ on the orbit space $M/G$. Our basic idea
here is to study the orbit space $M/G$ by examining holonomies of
loops passing through orbifold singularities of $M/G$. To be more
precise, we consider lifting a loop $\overline{\gamma}:S^1 @>>> M/G$
to a map $\gamma:\Bbb R @>>> M$, where $S^1=\Bbb R/\Bbb Z$. This lift
may not close after moving $1$ unit along $\Bbb R$ and the difference
between $\gamma(t)$ and $\gamma(t+1)$ comes from the action of an
element $g\in G$, the holonomy of $\overline{\gamma}$. When the loop
$\overline{\gamma}$ does not pass through orbifold points of $M/G$,
the conjugacy class of the holonomy is uniquely determined by
$\overline{\gamma}$. However, when the loop $\overline{\gamma}$ passes
through orbifold point, it can have lifts whose holonomies belong to
different conjugacy classes. Furthermore, it can have a lift whose
holonomy depend on the unit segment of $\Bbb R$ on which the holonomy
is measured. To avoid this complication, we consider only $g$-periodic
lifts. This is the notion of $g$-twisted free loop space $L_gM$
defined by
$$
L_gM=\{\gamma:\Bbb R @>>> M \mid \gamma(t+1)=g^{-1}\gamma(t),\
t\in\Bbb R\}.
\tag2-16
$$
Since any loop $\overline{\gamma}\in L(M/G)$ can be lifted to a
$g$-periodic lift for some $g\in G$, we have a surjective map 
$\coprod_{g\in G}L_gM @>>> L(M/G)$. On the space $\coprod_{g\in
G}L_gM$, the group $G$ acts inducing a homeomorphism $h\cdot:L_gM
@>{\cong}>> L_{hgh^{-1}}M$ for any $h,g\in G$. Quotienting by this
action, we get a surjective map
$$
\Bigl(\coprod_{g\in G}L_gM\Bigr)\!\!\Bigm/\!\!G\cong
\!\!\!\coprod_{[g]\in G_*}\!\!\!\bigl(L_gM/C_G(g)\bigr) 
@>{\text{onto}}>> L(M/G).
\tag2-17
$$
Here $G_*$ is the set of all conjugacy classes of $G$. This map is
$1:1$ on the subset of $L(M/G)$ consisting of loops not passing
through orbifold points. Thus, if the action of $G$ on $M$ is free,
then the above map is a homeomorphism. When a loop passes through
orbifold points, its inverse image is not unique, but finite,
corresponding to finitely many possibilities of different conjugacy
classes of lifts. Thus, in a sense, the above surjective map gives a
mild resolution of orbifold singularities.

Since $L_gM/C(g)$ is again an orbifold space, we can apply the above
procedure again on its free loop space. In fact, we can iterate this
procedure. To describe this general case, for any $d\ge1$ and for any
homomorphism $\phi:\Bbb Z^d @>>> G$, let $L_{\phi}M$ be the space of
twisted $d$-dimensional tori defined by 
$$
L_{\phi}M=\{\gamma:\Bbb R^d @>>> M \mid \gamma(\bold{t}+\bold{m})
=\phi(\bold{m})^{-1}\gamma(\bold{t}),\ \bold{t}\in
\Bbb R^d, \bold{m}\in\Bbb Z^d\ \}. 
\tag2-18
$$
Here $\phi$ plays a role of holonomy of the map $\overline{\gamma}:T^d
@>>> M/G$, where $T^d=\Bbb R^d/\Bbb Z^d$. Observe that any $\gamma$ in
$L_{\phi}M$ factors through the torus $T_{\phi}=\Bbb
R^d/\text{Ker}\,\phi$. As before, the action of any $h\in G$
induces a homeomorphism $h\cdot:L_{\phi}M @>{\cong}>> L_{h\phi
h^{-1}}M$. Let
$$
\Bbb {L}^d(M;G)=\Bigl(\!\!\!\!\!\!\!\!\!\!\!\!\!\!\!
\coprod_{\ \ \ \phi\in\text{Hom}(\Bbb Z^d,G)}
\!\!\!\!\!\!\!\!\!\!\!\!\!\!\!
L_{\phi}M\Bigr)\!\!\Bigm/\!\!G
=\!\!\!\!\!\!\!\!\!\!\!\!\!\!\!\!\!\!\!
\coprod_{\ \ \ \ [\phi]\in\text{Hom}(\Bbb Z^d,G)/G}
\!\!\!\!\!\!\!\!\!\!\!\!\!\!\!\!\!\!\!
\bigl(L_{\phi}M/C_G(\phi)\bigr).
\tag2-19
$$
Here, $C_G(\phi)\le G$ is the centralizer of $\phi$. We may call this
space $d$-th order twisted torus space for $(M;G)$. Let $\Bbb{T}=\Bbb
R^d/\cap_{\phi}\text{Ker}\,\phi$, where $\phi$ runs over all
homomorphisms $\text{Hom}(\Bbb Z^d,G)$. Then $\Bbb {T}$ is a
$d$-dimensional torus and it acts on $\Bbb{L}^d(M;G)$.

We have a canonical map $\Bbb{L}^d(M;G) @>>> L^d(M/G)$ from the above
space to $d$-th iterated free loop space on $M/G$. This map is no
longer surjective nor injective in general. Of course when the action
$G$ on $M$ is free, the above map is still a homeomorphism.  

The space $\Bbb{L}^d(M;G)$ can be thought of as the space of pairs
$(\overline{\gamma}, [\phi])$, where $\overline{\gamma}:T^d @>>> M/G$
is a $d$-torus in $M/G$, and $[\phi]$ is the conjugacy class of the
holonomy of a periodic lift of $\overline{\gamma}$ to a map
$\gamma:\Bbb R^d @>>> M$. 

We want to calculate ordinary Euler characteristic of the space
$\Bbb{L}^d(M;G)$. However, since this space is infinite dimensional,
it may have nonzero Betti numbers in arbitrarily high degrees. We
recall that for a finite dimensional manifold admitting a torus
action, it is well known that Euler characteristic of the fixed point
submanifold is the same as the Euler characteristic of the original
manifold. In fact, a formal application of Atiyah-Singer-Segal Fixed
Point Index Theorem predicts that the Euler characteristic of
$\Bbb{L}^d(M;G)$ must be the same as the Euler characteristic of
$\Bbb{T}$-fixed point subset. Thus, the Euler characteristic of
$\chi\bigl(\Bbb{L}^d(M;G)\bigr)$ ought to be given by
$$
\chi\bigl(\Bbb{L}^d(M;G)^{\Bbb{T}}\bigr)
=\!\!\!\!\!\!\!\!\!\!\!\!\!\!\!\!\!\!\!
\sum_{\ \ \ \ \ [\phi]\in\text{Hom}(\Bbb Z^d,G)/G}
\!\!\!\!\!\!\!\!\!\!\!\!\!\!\!\!\!\!\!
\chi\bigl(M^{\langle\phi\rangle}/C_G(\phi)\bigr)
=\chi^{(d)}(M;G). 
\tag2-20
$$
This is the geometric origin of our definition of higher order
orbifold Euler characteristic $\chi^{(d)}(M;G)$.

We can give a similar geometric interpretation of the higher order
$p$-primary orbifold Euler characteristic $\chi^{(d)}_p(M;G)$ as Euler
characteristic of an infinite dimensional twisted mapping space with a
torus action. The Euler characteristic of the fixed point subset under
this torus action is precisely given by $\chi^{(d)}_p(M;G)$, as in
(2-20). A similar consideration applies to generalized orbifold Euler
characteristic $\chi_K(M;G)$ for a general group $K$: we replace the
mapping space $\text{Map}(\Sigma, M/G)$, where $\Sigma$ is a manifold
with the fundamental group $\pi_1(\Sigma)=K$, by $G$-orbits of twisted
mapping spaces parametrized by $\text{Hom}(K,G)/G$, we then take the
Euler characteristic of constant maps.

\head
\S 3 Centralizers of wreath products
\endhead
This section reviews some facts on wreath products. In particular, we
explicitly describe the structure of centralizers of elements in
wreath products. This material may be well known to experts. For
example, the order of centralizers in wreath product is discussed in
Macdonald's book \cite{M2, p.171}. However, since the precise details
on this topic seem to be rather hard to locate in literature, our
explicit and direct description will make this paper more
self-contained and it may be useful for readers with different
expertise. The structure of centralizers is described in Theorem
3-5. This section is independent from the rest of the paper. 

Let $G$ be a finite group. The $n$-th symmetric group $\frak S_n$ acts
on the $n$-fold Cartesian product $G^n$ by $s(g_1,g_2,\dots,g_n)
=(g_{s^{-1}(1)},g_{s^{-1}(2)},\dots,g_{s^{-1}(n)})$, where $s\in\frak
S_n$, $g_i\in G$. The semidirect product defined by this action is the
wreath product $G\wr\frak S_n=G^n\rtimes \frak S_n$. We use the
notation $G_n$ to denote this wreath product. Product and inverse are
given by $(\bold{g}, s)(\bold{h},t)= (\bold{g}\cdot s(\bold{h}),st)$
and $(\bold{g},s)^{-1}=(s^{-1}(\bold{g}^{-1}), s^{-1})$.

\subhead
Conjugacy classes in wreath products
\endsubhead
Let $s=\prod_is_i$ be the cycle decomposition of $s$. If
$s=(i_1,i_2,\dots,i_r)$ is a linear representation of $s_i$, then the
product $g_{i_r}\cdots g_{i_2}g_{i_1}$ is called the cycle product of
$(\bold{g},s)$ corresponding to the above representation of the cycle
$s_i$. For each cycle $s_i$, the conjugacy class of its cycle product
is uniquely determined. 

Corresponding to $s_i$, let $\bold{g}_i$ be an element of $G^n$ whose
$a$-th component is given by $(\bold{g}_i)_a=(\bold{g})_a$ if
$a\in\{i_1,\cdots,i_r\}$, and $(\bold{g}_i)_a=1$ otherwise. Then for
each $i,j$, $(\bold{g}_i,s_i)$ and $(\bold{g}_j,s_j)$ commute and we
have $(\bold{g},s)=\prod_i(\bold{g}_i,s_i)$. We often write
$\bold{g}_i=(g_{i_1},g_{i_2},\dots,g_{i_r})$ as if it is an element of
$G^r$. The conjugacy class of the cycle product corresponding to $s_i$
is denoted by $[\bold{g}_i]$. 

Let $G_*$ denote the set of all conjugacy classes of $G$, and we fix
the representatives of conjugacy classes. For $[c]\in G_*$, let
$m_r(c)$ be the number of $r$-cycles in the cycle decomposition
$s=\prod_is_i$ whose cycle products belong to $[c]$. This yields a
partition-valued function $\rho:G_* @>>> \Cal P$, where $\Cal P$ is
the totality of partitions, defined by
$\rho([c])=(1^{m_1(c)}2^{m_2(c)}\cdots r^{m_r(c)}\cdots)$. Note that
$\sum_{[c],r}rm_r(c)=n$. The function $\rho$ associated to
$(\bold{g},s)\in G_n$ is called the type of $(\bold{g},s)$. It is well
known that the conjugacy class of $(\bold{g},s)$ in $G_n$ is
determined by its type. This can be explicitly seen using the
conjugation formulae in Proposition 3-1 below.

To describe details of the structure of the wreath product $G_n$, we
use the following notations. We express any element $(\bold{g},s)$ as
a product in two ways: 
$$
(\bold{g},s)=\prod_i(\bold{g}_i,s_i)
=\!\!\!\!\!\prod_{[c]\in G_*\ \ }\!\!\!\!\!\prod_{r\ge1}
\!\!\prod_{i=1}^{m_r(c)}
\!\!(\boldxi_{r,c,i},\sigma_{r,c,i}).
\tag3-1
$$
In the second expression, the conjugacy class of the cycle product
corresponding to $\sigma_{r,c,i}$ is $[\boldxi_{r,c,i}]=[c]$.

Suppose the conjugacy class of the cycle product
$[\bold{g}_i]=[g_{i_r}\cdots g_{i_1}]$ corresponding to $s_i$ is equal
to $[c]\in G_*$.  Choose and fix $p_i\in G$ such that $g_{i_r}\cdots
g_{i_1}=p_icp_i^{-1}$ for all $i$. Let
$$
\bolddelta(\bold{g}_i)
=(g_{i_1},g_{i_2}g_{i_1},\dots,g_{i_r}\cdots g_{i_1}),\quad
\Delta_i(p_i)=(p_i,\dots,p_i)\in G^r\!\!\overset{(i)}\to\subset \!\!G^n.
\tag3-2
$$
Here by $G^r\!\!\overset{(i)}\to\subset \!\!G^n$, we mean that the components
of the above elements at the position $a\notin\{i_1,i_2,\dots,i_r\}$
are $1\in G$. This convention applies throughout this section. In the
above, $\Delta_i$ is the diagonal map along the components of
$s_i$. Let $\bold{c}_{r}=(c,1,\dots,1)\in G^r$ and $\bold{c}_{r,i}\in
G^n$ is $\bold{c}_r$ along the positions which appear as components of
$\sigma_{r,c,i}$ or $s_i$. The following proposition can be checked by
direct calculation.

\proclaim{Proposition 3-1} $(1)$ Let
$(\bold{g},s)=\prod_i(\bold{g}_i,s_i)$. For a given $i$, suppose
$[\bold{g}_i]=[c]$ and $|s_i|=r$. Then $\bolddelta_i
=\bolddelta(\bold{g}_i)\cdot\Delta_i(p_i)\in 
G^r\!\!\overset{(i)}\to\subset \!\!G^n$ has the property 
$$
(\bold{g}_i,s_i)=(\bolddelta_i,1)\bigl((c,1,\cdots,1),s_i\bigr)
(\bolddelta_i, 1)^{-1}. 
\tag3-3
$$
For $i\not=i$, elements $(\bolddelta_i,1)$ and
$(\bolddelta_j,1)$ commute. 

$(2)$ Let
$(\bolddelta,1)=\prod_i(\bolddelta_i,1)$. Then 
$$
(\bold{g},s)=(\bolddelta,1)\Bigl(\prod_{[c]}\prod_{r\ge1}
\!\!\prod_{i=1}^{m_r(c)}\!\!(\bold{c}_{r,i},\sigma_{r,c,i})\Bigr)
(\bolddelta,1)^{-1}.
\tag3-4
$$
\endproclaim

>From this, it is clear that the conjugacy class of $(\bold{g},s)\in G^n$
is determined by its type $\{m_r(c)\}_{r,[c]}$. 

\subhead
Actions of wreath products
\endsubhead
Let $M$ be a $G$-manifold. The wreath product $G_n$ acts on the
$n$-fold Cartesian product $M^n$ by
$$
\bigl((g_1,g_2,\dots,g_n),s\bigr)(x_1,x_2,\dots,x_n)
=(g_1x_{s^{-1}(1)},g_2x_{s^{-1}(2)},\dots,g_nx_{s^{-1}(n)}).
$$ 
The above conjugation formula shows that the fixed point subset
of $M^n$ under the action of $(\bold{g},s)$ is completely determined
by its type. 

\proclaim{Proposition 3-2} Suppose an element $(\bold{g},s)\in G_n$ is
of type $\{m_r(c)\}$. Then
$$
(M^n)^{\langle(\bold{g},s)\rangle}\cong\prod_{[c]}
\bigl(M^{\langle c\rangle}\bigr)^{\sum_rm_r(c)}.
\tag3-5
$$
\endproclaim
\demo{Proof} With respect to the decomposition (3-1) of
$(\bold{g},s)$, we have $(M^n)^{\langle(\bold{g},s)\rangle}\cong
\prod_i(M^{r_i})^{\langle(\bold{g}_i,s_i)\rangle}$, where
$r_i=|s_i|$, and $M^{r_i}\subset M^n$ corresponds to components of
$s_i$. If $[\bold{g}_i]=[c]$, then we have the following isomorphisms
$$
M^{\langle c\rangle} @>{\cong}>{\Delta_i}> 
(M^{r_i})^{\langle(c_{r_i},s_i)\rangle}
@>{\cong}>{(\bolddelta_i,1)}> 
(M^{r_i})^{\langle(\bold{g}_i,s_i)\rangle}.
$$
Now in terms of the other decomposition of $(\bold{g},s)$ in (3-1),
since $[\boldxi_{r,c,i}]=[c]$, the above isomorphism implies
$$
(M^n)^{\langle(\bold{g},s)\rangle}\cong
\prod_{[c]}\prod_{r\ge1}\!\!\prod_{i=1}^{m_r(c)}\!\!
(M^{r_i})^{\langle(\boldxi_{r,c,i},\sigma_{r,c,i})\rangle}
\cong\prod_{[c]}(M^{\langle c\rangle})^{\sum_rm_r(c)}.
$$
This completes the proof. 
\qed
\enddemo

\subhead
Centralizers in wreath products
\endsubhead
Next, we describe the centralizer $C_{G_n}\!\bigl((\bold{g},s)\bigr)$
in the wreath product $G_n$. Let $(\bold{h},t)\in
C_{G_n}\!\bigl((\bold{g},s)\bigr)$. Then
$(\bold{h},t)(\bold{g},s)(\bold{h},t)^{-1}=(\bold{g},s)$. In terms of
the cycle decomposition $(\bold{g},s)=\prod_i(\bold{g}_i,s_i)$, we see
that for each $i$ there exits a unique $j$ such that
$(\bold{h},t)(\bold{g}_i,s_i)
(\bold{h},t)^{-1}=(\bold{g}_j,s_j)$. Since the conjugation preserves
the type, we must have $[\bold{g}_i]=[\bold{g}_j]=[c]$ for some
$[c]\in G_*$, and $|s_i|=|s_j|$. Thus with respect to the second
decomposition in (3-1), the conjugation by $(\bold{h},t)$ permutes
$m_r(c)$ elements $\{(\boldxi_{r,c,i},\sigma_{r,c,i})\}_i$
for each $[c]\in G_*$ and $r\ge1$. Thus we have a homomorphism
$$
p: C_{G_n}\!\bigl((\bold{g},s)\bigr) @>{\hphantom{longarr}}>> 
\prod_{[c]}\prod_{r\ge1}\frak S_{m_r(c)}.
\tag3-6
$$

\proclaim{Lemma 3-3} The above homomorphism $p$ is split
surjective. 
\endproclaim
\demo{Proof} First we construct a homomorphism
$\lambda: \prod_{[c]}\prod_{r}\frak S_{m_r(c)} @>>> \frak S_n$ as
follows. In the decomposition $(\bold{g},s)=\prod_i(\bold{g}_i,s_i)$,
we write each cycle $s_i$ starting with the smallest integer. If
$\overline{t}\in\prod_{[c]}\prod_{r}\frak S_{m_r(c)}$ sends the cycle
$s_i=(i_1,i_2,\dots,i_r)$ to $s_j=(j_1,j_2,\dots,j_r)$, then let
$\lambda(\overline{t})=t\in\frak S_n$ where $t(i_{\ell})=j_{\ell}$,
$1\le\ell\le r$. It is clear that $\lambda$ defines a homomorphism,
and any element in the image of $\lambda$ commutes with $s$, because
$t$ permutes cycles preserving the smallest integers in cycles
appearing in the cycle decomposition of $s$. Previously, we
constructed an element $\bolddelta\in G_n$ for each element
$(\bold{g},s)$ in Proposition 3-1. Define a homomorphism $\Lambda:
\prod_{[c]}\prod_{r}\frak S_{m_r(c)} @>>>
C_{G_n}\!\bigl((\bold{g},s)\bigr)$ by $\Lambda(\overline{t})
=(\bolddelta,1)\bigl(1,\lambda(\overline{t})\bigr)
(\bolddelta,1)^{-1}$. Using (3-4), we can check that
$\Lambda(\overline{t})$ commutes with $(\bold{g},s)$ for any
$\overline{t}$. Since conjugation by $\Lambda(\overline{t})$ on
$(\bold{g},s)$ induces the permutation $\overline{t}$ among
$\{(\bold{g}_i,s_i)\}_i$, we have $p\circ\Lambda=\text{identity}$ and
$\Lambda$ is a splitting of $p$. This completes the proof.
\qed
\enddemo

Next we examine the kernel of the homomorphism $p$ in (3-6). If
$(\bold{h},t)\in\text{Ker}\,p$, then
$(\bold{h},t)(\bold{g}_i,s_i)(\bold{h},t)^{-1}=(\bold{g}_i,s_i)$ for
all $i$. In particular, we have $ts_it^{-1}=s_i$ for all $i$, and
consequently $t$ must be a product of powers of $s_i$'s. Thus, we may
write
$(\bold{h},t)=\prod_i\bigl\{(\bold{h}_i,1)(\bold{g}_i,s_i)^{k_i}\bigr\}$
for some $\bold{h}_i\in G^{|s_i|}\subset G^n$ and $0\le k_i< |s_i|$,
where $(\bold{h}_i,1)$ commutes with $(\bold{g}_i,s_i)$ for any
$i$. Let $G_r^{(i)}=G^r\rtimes\frak S_r\subset G_n$ be a subgroup of
$G_n$ isomorphic to $G_r$ corresponding to positions appearing in
$s_i$. Recall that we defined $\bolddelta_i$ in (3-3).

\proclaim{Lemma 3-4} For a given $i$, suppose the cycle product
corresponding to $s_i$ is such that $[\bold{g}_i]=[c]$. Then
$$
\aligned
C_{G_r^{(i)}}\bigl((\bold{g}_i,s_i)\bigr)
&=\{(\bold{h}_i,1)\cdot(\bold{g}_i,s_i)^{k_i}\mid 
0\le k_i<|s_i|, [(\bold{h}_i,1),(\bold{g}_i,s_i)]=1\}\\
&=\bigl(\bolddelta_i\cdot\Delta_i(C_G(c))
\cdot\bolddelta_i^{-1},1\bigr)\cdot
\bigl\langle(\bold{g}_i,s_i)\bigr\rangle\\
&\cong C_G(c)\cdot\langle a_{r,c}\rangle,\ \ \text{where } (a_{r,c})^r=c
\text{ and }[a_{r,c},C_G(c)]=1.
\endaligned
\tag3-7
$$
Here $a_{r,c}=\bigl(\bold{c}_r,(12\cdots r)\bigr)$. 
\endproclaim
\demo{Proof} The first equality is obvious. For the second one, first
observe that if $(\bold{h},1)$ with $\bold{h}\in G^r$ commutes with
$\bigl(\bold{c}_r,(12\cdots r)\bigr)\in G_r$, then $\bold{h}$ must be of the
form $(h,h,\dots,h)\in G^r$ with $h\in C_G(c)$. Conjugation by
$\bolddelta_i$ gives the second description. For the third
description, we simply observe that $\bigl(\bold{c}_r,(12\cdots
r)\bigr)^r=\bigl(\Delta(h),1\bigr)\in G_r$. This completes the proof.
\qed
\enddemo
Thus we have a split exact sequence 
$$
1 @>>> \prod_{[c]}\prod_{r\ge1}\!\!\prod_{i=1}^{m_r(c)}\!\!
C_{G_r^{(i)}}
\bigl((\boldxi_{r,c,i},\sigma_{r,c,i})\bigr) @>>> 
C_{G_n}\!\bigl((\bold{g},s)\bigr) @>{p}>> 
\prod_{[c]}\prod_{r\ge1}\frak S_{m_r(c)} @>>> 1.
\tag3-8
$$
Since conjugation preserves the type of elements in $G_n$, the
centralizer splits into $(r,[c])$-components
$$
C_{G_n}\!\bigl((\bold{g},s)\bigr)=\prod_{[c]}\prod_{r\ge1}
C_{G_n}\!\bigl((\bold{g},s)\bigr)_{(r,[c])},
$$
where $C_{G_n}\!\bigl((\bold{g},s)\bigr)_{(r,[c])}$ is the centralizer
of the element $\prod_i(\boldxi_{r,c,i},\sigma_{r,c,i})$ in
the subgroup $G_{rm_r(c)}$ corresponding to positions appearing in
$\sigma_{r,c,i}$ for $1\le i\le m_r(c)$. The conjugation by
$(\bolddelta_i,1)$ maps the following split exact sequence
$$
1 @>>> \!\!\prod_{i=1}^{m_r(c)}\!\!C_{G_r^{(i)}}
\bigl((\boldxi_{r,c,i},\sigma_{r,c,i})\bigr) @>>> 
C_{G_n}\!\bigl((\bold{g},s)\bigr)_{(r,[c])} @>{p_{r,c}}>> 
\frak S_{m_r(c)} @>>> 1 
\tag3-9
$$
isomorphically into the following split exact sequence 
$$
1 @>>> \!\!\prod_{i=1}^{m_r(c)}\!\!\Delta_i\bigl(C_G(c)\bigr)\cdot
\langle a_{r,c}^{(i)}\rangle @>>> 
C_{G_{rm_r(c)}}\Bigl(\prod_{i=1}^{m_r(c)}
(\bold{c}_{r,i},\sigma_{r,c,i})\Bigr)
@>>> \frak S_{m_r(c)} @>>> 1.
$$
Here $a_{r,c}^{(i)}$ is $a_{r,c}$ along components of
$\sigma_{r,c,i}$. Direct calculation shows that in the second exact
sequence, $\frak S_{m_r(c)}$ acts on the left side product by
permuting factors. Hence the semidirect product structure in (3-9) is
indeed isomorphic to a wreath product. Hence we obtain the following
description of the centralizer of $(\bold{g},s)$ in $G_n$.

\proclaim{Theorem 3-5} Let $(\bold{g},s)\in G_n$ have type
$\{m_r(c)\}_{r,[c]}$. Then
$$
C_{G_n}\!\bigl((\bold{g},s)\bigr)\cong\prod_{[c]}\prod_{r\ge1}
\bigl\{\bigl(C_G(c)\cdot\langle a_{r,c}\rangle\bigr)
\wr\frak S_{m_r(c)}\bigr\},
\tag3-10
$$
where $(a_{r,c})^r=c\in C_G(c)$. Here, the isomorphism is induced by
conjugation by $\bolddelta$ in Proposition {\rm 3-1}.
\endproclaim

\head
\S 4 Higher order orbifold Euler characteristic of symmetric
products
\endhead

In this section, we prove Theorem A in the introduction. Explicitly
writing out $j_r(\Bbb Z^d)$, the formula we prove is the following:
$$
\sum_{n\ge0}q^n\chi^{(d)}(M;G\wr\frak S_n)
=\Bigl[\!\!\!\!\!\!\!\!
\prod_{r_1,r_2,\dots,r_d\ge1}\!\!\!\!\!\!\!\!
(1-q^{r_1r_2\cdots r_d})
^{r_2r_3^2\cdots r_d^{d-1}}\Bigr]^{(-1)\chi^{(d)}(M;G)}.
\tag4-1
$$
When $d=0$, this is Macdonald's formula applied to $M/G$: 
$$
\sum_{n\ge0}q^n\chi\bigl(SP^n(M/G)\bigr)
=\frac1{(1-q)^{\chi(M/G)}}.
\tag4-2
$$
When $d=1$, the formula was proved by Wang \cite{W}. We prove formula
(4-1) by induction on $d\ge0$, using Macdonald's formula as the start
of induction. For the inductive step, we need the following
Lemma. 

\proclaim{Lemma 4-1}  Let $G\cdot\langle a\rangle$ be a group generated
by a finite group $G$ and an element $a$ such that $a$ commutes with
any element of $G$ and $\langle a\rangle\cap G=\langle a^r\rangle$ for
some integer $r\ge1$. Suppose the element $a$ acts trivially on a
$G$-manifold $M$. Then
$$
\chi^{(d)}(M;G\cdot\langle a\rangle)=r^d\chi^{(d)}(M;G).
\tag4-3
$$
\endproclaim
\demo{Proof} First note that $G\cdot\langle
a\rangle=\coprod_{i=0}^{r-1}G\cdot a^i$ and $|G\cdot\langle a\rangle|
=r\cdot|G|$. Observe that two elements of the form $ga^i$ and $ha^j$,
where $g,h\in G$, commute if and only if $g,h$ commute, since $a$ is
in the center of $G\cdot\langle a\rangle$. Now by definition, 
$$
\chi^{(d)}(M;G\cdot\langle a\rangle)
=\frac1{r\cdot|G|}\!\!\!\!\!\!\!\!\!\!\!\!\!
\sum \Sb \ \ \ \ \ \ (g_1,\dots,g_{d+1}) \\ \ \ \ 0\le i_{\ell}<r \endSb
\!\!\!\!\!\!\!\!\!\!\!\!\!
\chi\bigl(M^{\langle g_1a^{i_1},g_2a^{i_2},\dots,g_{d+1}a^{i_{d+1}}
\rangle}\bigr),
$$
where $(g_1,\dots,g_{d+1})$ runs over all $(d+1)$-tuples of mutually
commuting elements of $G$, and the index $\ell$ runs over $1\le
\ell\le d+1$. Since the element $a$ acts trivially on $M$, the fixed
point subset above is the same as $M^{\langle
g_1,\dots,g_{d+1}\rangle}$. Hence summing over $i_{\ell}$'s first, the
above becomes
$$
\chi^{(d)}(M;G\cdot\langle a\rangle)=\frac{r^{d+1}}{r\cdot|G|}
\!\!\!\!\!\!\!\!\!\!\!\!\!\!\!
\sum_{\ \ \ \ \ \ \ (g_1,\dots,g_{d+1})}
\!\!\!\!\!\!\!\!\!\!\!\!\!\!\!
\chi\bigl(M^{\langle g_1,\dots,g_{d+1}\rangle}\bigr)
=r^d\cdot\chi^{(d)}(M;G). 
$$
This completes the proof. 
\qed
\enddemo

\demo{Proof of formula {\rm (4-1)}} By induction on $d\ge0$. When
$d=0$, the formula is Macdonald's formula (4-2) and hence it is
valid. 

Assume the formula is valid for $\chi^{(d-1)}$ for $d\ge1$. Let
$G_*=\{[c]\}$ be the totality of conjugacy classes of $G$. By
Proposition 2-5, we have 
$$
\sum_{n\ge0}q^n\chi^{(d)}(M;G_n)
=\sum_{n\ge0}q^n\!\!\sum_{[(\bold{g},s)]}\!\!
\chi^{(d-1)}\bigl((M^n)^{\langle(\bold{g},s)\rangle};
C_{G_n}\!((\bold{g},s))\bigr).
\tag{$*$}
$$
Let $(\bold{g},s)\in G_n$ have type $\{m_r(c)\}$. Then by Proposition
3-2 and Theorem 3-5, we have the following compatible isomorphisms:
$$
\align
(M^n)^{\langle(\bold{g},s)\rangle}&\cong \prod_{[c]}\prod_{r\ge1}
(M^{\langle c\rangle})^{m_r(c)},\\
C_{G_n}\!\bigl((\bold{g},s)\bigr)&\cong\prod_{[c]}\prod_{r\ge1}
\bigl\{(C_G(c)\cdot\langle a_{r,c}\rangle)\wr\frak S_{m_r(c)}\bigr\},
\endalign
$$
where $(a_{r,c})^r=c\in C_G(c)$ and $a_{r,c}$ acts trivially on
$M^{\langle c\rangle}$. The above isomorphisms are compatible in the
sense that the action of $C_{G_n}\!\bigl((\bold{g},s)\bigr)$ on 
$(M^n)^{\langle(\bold{g},s)\rangle}$ translates, via conjugation by
$(\bolddelta,1)$ in Proposition 3-1, to the action of 
the wreath product $(C_G(c)\cdot\langle a_{r,c}\rangle)\wr\frak
S_{m_r(c)}$ on $(M^{\langle c\rangle})^{m_r(c)}$ for any $[c]\in G_*$
and $r\ge1$. Since the conjugacy classes of elements in $G_n$ are
determined by their types, the summation over all conjugacy classes 
$[(\bold{g},s)]$ corresponds to the summation over all $m_r(c)\ge0$ for
all $[c]\in G_*$ and $r\ge1$ subject to $\sum_{[c],r}rm_r(c)=n$. By
the multiplicativity of generalized orbifold Euler characteristic
(2-3), the formula $(*)$ becomes 
$$
\align
\sum_{n\ge0}q^n\chi^{(d)}&(M^n;G_n)\\
&=\sum_{n\ge0}q^n\!\!\!\!\!\!\!\!\!
\sum \Sb m_r(c)\ge0 \\\sum rm_r(c)=n \endSb
\!\!\!\prod_{[c],r}\chi^{(d-1)}\bigl((M^{\langle c\rangle})^{m_r(c)};
(C_G(c)\langle a_{r,c}\rangle)\wr\frak S_{m_r(c)}\bigr)\\
&=\!\!\!\!\!\!\!\!\sum_{m_r(c)\ge0\ \ \ \ }\!\!\!\!\!\!\!\!
\prod_{[c],r}(q^r)^{m_r(c)}\chi^{(d-1)}\bigl(
(M^{\langle c\rangle})^{m_r(c)};
(C_G(c)\langle a_{r,c}\rangle)\wr\frak S_{m_r(c)}\bigr)\\
&=\prod_{[c],r}\sum_{m\ge0}(q^r)^m\chi^{(d-1)}
\bigl((M^{\langle c\rangle})^{m};
(C_G(c)\cdot\langle a_{r,c}\rangle)_m\bigr)
\intertext{By inductive hypothesis, this is equal to }
&=\prod_{[c],r}\Bigl[\!\!\!\!\!\!\!\!\!\!\!\!\!\!\!
\prod_{\ \ \ \ \ \ \ \  r_1,\dots,r_{d-1}\ge1}
\!\!\!\!\!\!\!\!\!\!\!\!\!\!\!\!\!\!
\bigl(1-(q^r)^{r_1\dots r_{d-1}}
\bigr)^{r_2r_3^2\cdots r_{d-1}^{d-2}}\Bigr]
^{-\chi^{(d-1)}(M^{\langle c\rangle};C_G(c)\langle a_{r,c}\rangle)}\\
&=\Bigl[\!\!\!\!\!\!\!\!\!\prod_{r,r_1,\dots,r_{d-1}}
\!\!\!\!\!\!\!\!
\bigl(1-q^{rr_1\cdots r_{d-1}}\bigr)^{r_2r_3^2\cdots r_{d-1}^{d-2}}
\Bigr]^{-\sum_{[c]}\chi^{(d-1)}(M^{\langle
c\rangle};C_G(c)\langle a_{r,c}\rangle)}.
\endalign
$$
By Lemma 4-1, $\chi^{(d-1)}(M^{\langle c\rangle};
C_G(c)\cdot\langle a_{r,c}\rangle)=r^{d-1}\chi^{(d-1)}
(M^{\langle c\rangle};C_G(c))$. Hence summing over $[c]\in G_*$ 
and using Proposition 2-5, we see that the exponent is equal to 
$(-1)r^{d-1}\chi^{(d)}(M;G)$. Thus, renaming $r$ as $r_d$, the above
is equal to the right hand side of (4-1). This completes the inductive
step and the proof is complete. 
\qed
\enddemo

Now let $M=\text{pt}$. Using (2-2) with $K=\Bbb Z^d$ and $G$
replaced by $G$ or $G_n$, we get 

\proclaim{Corollary 4-2}  For each $d\ge0$ and for any finite group
$G$, we have 
$$
\sum_{n\ge0}q^n\bigl|\text{\rm Hom}\bigl(\Bbb Z^d, G_n\bigr)/G_n\bigr|
=\Bigl[\!\!\!\!\!\!\prod_{r_1,\dots,r_d\ge1}\!\!\!\!\!\!
(1-q^{r_1\cdots r_d})^{r_2r_3^2\cdots r_d^{d-1}}\Bigr]^{(-1)
|\text{Hom}(\Bbb Z^d,G)/G|}.
\tag4-4
$$
\endproclaim

The above formula is the formula (1-7) in the
introduction. Furthermore, letting $G$ be the trivial group, we get 
$$
\sum_{n\ge0}q^n\bigl|\text{\rm Hom}
\bigl(\Bbb Z^d,\frak S_n\bigr)/\frak S_n\bigr|
=\Bigl[\!\!\!\!\!\!\prod_{r_1,\dots,r_d\ge1}\!\!\!\!\!\!
(1-q^{r_1\cdots r_d})
^{r_2r_3^2\cdots r_d^{d-1}}\Bigr]^{(-1)}.
\tag4-5
$$
Here, as remarked in the introduction, we recognize $|\text{\rm Hom}
\bigl(\Bbb Z^d,\frak S_n\bigr)/\frak S_n\bigr|$ as the number of
isomorphism classes of $\Bbb Z^d$-sets of order $n$. Any finite $\Bbb
Z^d$-set decomposes into a union of transitive $\Bbb Z^d$-sets, and
any isomorphism class of transitive $\Bbb Z^d$-set of order $r$
corresponds to a unique subgroup of $\Bbb Z^d$ of index $r$, by taking
the isotropy subgroup. Thus, letting $j_r(\Bbb Z^d)$ be the number of
index $r$ subgroups of $\Bbb Z^d$, we have
$$
\sum_{n\ge0}q^n\bigl|\text{\rm Hom}
\bigl(\Bbb Z^d,\frak S_n\bigr)/\frak S_n\bigr|
=\Bigl[\prod_{r\ge1}(1-q^r)^{j_r(\Bbb Z^d)}\Bigr]^{(-1)}.
\tag4-6
$$
By comparing (4-5) and (4-6), we get a formula for $j_r(\Bbb
Z^d)$. However, we can easily directly calculate the number $j_r(\Bbb
Z^d)$ as follows. This calculation is well known (for its history, see
\cite{So}) and gives an alternate proof of (4-5).

\proclaim{Lemma 4-4} For any $r\ge1$, $d\ge1$, we have 
$$
j_r(\Bbb Z^d)=\!\!\!\!\!\!\sum_{r_1\cdots r_d=r}
\!\!\!\!\!\!r_2r_3^2\cdots r_d^{d-1},\quad\text{and}\quad
j_r(\Bbb Z^d)=\sum_{m|r}m\cdot j_m(\Bbb Z^{d-1}). 
\tag4-7
$$
\endproclaim
\demo{Proof} Let $e_1,e_2,\dots,e_d$ be the standard basis of the
lattice. It is easy to see that any sublattice of index $r$ has a
unique basis $\{x_i\}_{i=1}^d$ of the form $x_i=r_ie_i+\sum_{i<j\le
d}a_{ij}e_j$ for $1\le i\le d$, where $r_1r_2\cdots r_d=r$ and $0\le
a_{ij}<r_j$ for $1\le i<j\le d$. For any choice of $\{a_{ij}\}$
satisfying the condition, there exists a sublattice of rank $d$ of
index $r=r_1\cdots r_d$. Since given $r_1,\dots,r_d\ge1$, there are
$r_2r_3^2\cdots r_d^{d-1}$ choices of $a_{ij}$'s, the total number of
sublattices of index $r$ in $\Bbb Z^d$ is given by $\sum_{r_1\cdots
r_d=r}r_2r_3^2\cdots r_d^{d-1}$. 

The second equality easily follows from the first. This completes the
proof.
\qed
\enddemo

The proof of the identity (1-10) mentioned in the introduction can be
readily proved using (4-7) by induction on $d\ge1$. We will discuss
the corresponding $p$-local situation in the next section.

\head
\S 5 Euler characteristic of equivariant Morava K-theory of symmetric
products
\endhead

As before, let $G$ be a finite group and let $M$ be a
$G$-manifold. Let $d\ge0$ be an integer. The equivariant $d$-th Morava
$K$-theory of $M$ is defined as $K(d)^*_G(M)=K(d)^*(EG\times_G
M)$. Since $K(d)^*=\Bbb F_p[v_d,v_d^{-1}]$ with
$|v_d|=-2(p^d-1)$ is a graded field, any $K(d)^*$-module is free and
we can count the dimension over $K(d)^*$, if it is finite. Let
$\chi\bigl(K(d)^*_G(M)\bigr)$ be its Euler characteristic defined by
$$
\chi\bigl(K(d)^*_G(M)\bigr)=
\dim K(d)^{\text{even}}(EG\times_G M)-
\dim K(d)^{\text{odd}}(EG\times_G M).
\tag5-1
$$
When $d=0$, we have $K(0)^*(\ \cdot\ )=H^*(\ \cdot\ ;\Bbb Q)$ and
$\chi\bigl(K(0)^*_G(M)\bigr)=\chi(M/G)$, the ordinary Euler
characteristic of the orbit space. Hopkins, Kuhn, and Ravenel compute
this number (5-1) using a general theory of complex oriented additive
functions \cite{HKR, Theorem 4.12}. Their result is
$$
\chi\bigl(K(d)^*_G(M)\bigr)=\frac1{|G|}
\!\!\!\!\!\sum \Sb A\le G \\ A:\text{abelian} \endSb 
\!\!\!\!\!
|A|\cdot|A_{(p)}|^d\cdot \mu_A^{\Bbb C}(M),
\tag5-2
$$
where the M\"obius function $\mu_A^{\Bbb C}(M)$ is defined in (2-7),
and the summation is over all abelian subgroups of $G$. Now comparing
(5-2) with our calculation of higher order $p$-primary orbifold Euler
characteristic of $(M;G)$ in Proposition 2-7, we realize that these
two quantities are in fact equal for all $d\ge0$.

\proclaim{Proposition 5-1} For any $G$-manifold $M$, and for any
$d\ge0$, 
$$
\chi\bigl(K(d)^*_G(M)\bigr)=\chi^{(d)}_p(M;G)
=\frac1{|G|}\!\!\!\!\!\!\!\!\!\!\!\!
\sum_{\ \ \ \ \ \phi:\Bbb Z_p^d\times\Bbb Z @>>> G}
\!\!\!\!\!\!\!\!\!\!\!\!\!\!\chi(M^{\langle\phi\rangle}). 
\tag5-3
$$
\endproclaim

Our objective in this section is to calculate the Euler characteristic
of equivariant Morava K-theory of symmetric products $(M^n;G\wr\frak
S_n)$ for $n\ge1$. By Proposition 5-1, this homotopy theoretic number
can be calculated as the higher order $p$-primary orbifold Euler
characteristic. We prove 

\proclaim{Theorem 5-2} For any $d\ge0$ and $G$-manifold $M$, 
$$
\sum_{n\ge0}q^n \chi^{(d)}_p(M^n;G\wr\frak S_n)
=\Bigl[\!\!\!\!\!\!\!\prod_{\ell_1,\dots,\ell_d\ge0}
\!\!\!\!\!\!\!
(1-q^{p^{\ell_1}p^{\ell_2}\cdots p^{\ell_d}})
^{p^{\ell_2}p^{2\ell_3}\cdots p^{(d-1)\ell_d}}\Bigr]
^{(-1)\chi^{(d)}_p(M;G)}.
\tag5-4
$$
\endproclaim
The proof is very similar to the one for the formula (4-1). Note the
similarities of exponents. Since the formula (5-4) is $p$-primary,
there are differences at many parts of the proof, although the idea of
the proof is the same. Thus, we believe that it is better to give a
complete proof of the above formula (5-4) rather than explaining
differences of proofs between formulae (4-1) and (5-4).

\demo{Proof of Theorem {\rm 5-2}} By induction on $d\ge0$. When
$d=0$, by (2-13) we have $\chi^{(0)}_p(M^n;G\wr\frak S_n)=
\chi\bigl(SP^n(M/G)\bigr)$, and the formula (5-4) in this case asserts
$$
\sum_{n\ge0}q^n\chi\bigl(SP^n(M/G)\bigr)
=\frac1{(1-q)^{\chi(M/G)}},
$$
which is valid due to Macdonald's formula (1-8).

Assume that the formula (5-4) is valid for $\chi^{(d-1)}_p$ for $d\ge1$. By
Proposition 2-6, 
$$
\sum_{n\ge0}q^n\chi^{(d)}_p(M^n;G\wr\frak S_n)
=\sum_{n\ge0}q^n\sum_{[\alpha]}\chi^{(d-1)}_p
\bigl(M^{\langle\alpha\rangle};C_{G_n}\!(\alpha)\bigr).
\tag{$*$}
$$
Here the second summation on the right hand side is over all
$G_n$-conjugacy classes $[\alpha]$ of elements of $p$-power order in
$G_n$, that is $[\alpha]\in\text{Hom}(\Bbb Z_p,G_n)/G_n$. Let
$\alpha=(\bold{g},s)\in G_n$. Since $\alpha$ has order a power of $p$,
the second component $s\in\frak S_n$ must have order a power of
$p$. Thus, the type of $\alpha=(\bold{g},s)$ must be of the form
$\{m_{p^r}(c)\}_{r,[c]}$. Here, $[c]$ runs over all $G$-conjugacy
classes of elements of order powers of $p$. We indicate this by the
notation $[c]_p$. So $[c]_p\in\text{Hom}(\Bbb Z_p,G)/G$. By
Proposition 3-1, $\alpha$ is conjugate to an element of the form
$$
\alpha\sim \prod_{r\ge0}\prod_{[c]_p}\!\!\!\!\prod_{i=1}^{m_{p^r}(c)\ }
\!\!\!\!\oversetbrace{a_{p^r,c,i}}\to{\bigl(\undersetbrace{p^r}
\to{(c,1,\dots,1)},
\sigma_{p^r,c,i})}, \quad\text{where } (a_{p^r,c,i})^{p^r}
\!\!\!=\bigl(\undersetbrace{p^r}\to{(c,c,\dots,c)},1\bigr).
$$
By Proposition 3-2 and Theorem 3-5, the fixed point subset of $M^n$
under the action of $\alpha$, and the centralizer of $\alpha$ in $G_n$
are each isomorphic to
$$
\align
(M^n)^{\langle\alpha\rangle}&\cong\prod_{[c]_p}
(M^{\langle c\rangle})^{\sum_{r}m_{p^r}(c)},\\
C_{G_n}\!(\alpha)\cong\prod_{[c]_p}\prod_{r\ge0}&
\bigl\{(C_G(c)\cdot\langle a_{p^r,c}\rangle)\wr
\frak S_{m_{p^r}(c)}\bigr\},\quad (a_{p^r,c})^{p^r}\!\!=c\in C_G(c). 
\endalign
$$
The above isomorphisms are compatible with the action of the
centralizer on the fixed point subset. The summation over all
conjugacy classes $[\alpha]_p$ can be replaced by the summation over
all the types $\{m_{p^r}(c)\}_{r,[c]_p}$. By multiplicativity of
$\chi^{(d-1)}_p$, the right hand side of $(*)$ becomes 
$$
\align
(*)&=\sum_{n\ge0}q^n\!\!\!\!\!\!\!\!
\sum \Sb m_{p^r}(c)\ge0 \\ \sum p^rm_{p^r}(c)=n \endSb
\!\!\!\!
\prod_{[c]_p,r}\chi^{(d-1)}_p\bigl(
(M^{\langle c\rangle})^{m_{p^r}(c)};(C_G(c)\cdot
\langle a_{p^r,c}\rangle)\wr\frak S_{m_{p^r}(c)}\bigr)\hbox to 23pt{\hss}\\
&=\!\!\!\!\!\!\!\!\sum_{m_{p^r}(c)\ge0\ \ \ }\!\!\!\!\!\!\!
\prod_{\ [c]_p,r}\!\!(q^{p^r})^{m_{p^r}(c)}
\chi^{(d-1)}_p\bigl(
(M^{\langle c\rangle})^{m_{p^r}(c)};(C_G(c)\cdot
\langle a_{p^r,c}\rangle)\wr\frak S_{m_{p^r}(c)}\bigr)\\
&=\prod_{[c]_p,r}\sum_{m\ge0}(q^{p^r})^m
\chi^{(d-1)}_p\bigl(
(M^{\langle c\rangle})^{m};(C_G(c)\cdot
\langle a_{p^r,c}\rangle)\wr\frak S_{m}\bigr)
\endalign$$
By inductive hypothesis, the summation inside is given by 
$$\align
&=\!\!\prod_{[c]_p,r}\Bigl[\!\!\!\!\!\!\!\!\!\!\!\!\!\!
\prod_{\ \ \ \ \ \ \ell_1,\dots,\ell_{d-1}\ge0}\!\!\!\!\!\!\!\!\!\!\!\!\!
(1-(q^{p^r})^{p^{\ell_1}\cdots p^{\ell_{d-1}}})^{p^{\ell_2}p^{2\ell_3}
\cdots p^{(d-2)\ell_{d-1}}}\Bigr]
^{(-1)\chi^{(d-1)}_p(M^{\langle c\rangle};C_G(c)\cdot
\langle a_{p^r,c}\rangle)}\\
&=\Bigl[\!\!\!\!\!\!\!\!\!\!\prod_{\ell_1,\dots,\ell_{d-1},r\ge0}
\!\!\!\!\!\!\!\!\!\!\!
(1-q^{p^{\ell_1}\cdots p^{\ell_{d-1}}p^r})^{p^{\ell_2}p^{2\ell_3}
\cdots p^{(d-2)\ell_{d-1}}}\Bigr]
^{(-1)\sum_{[c]_p}\chi^{(d-1)}_p(M^{\langle c\rangle};C_G(c)\cdot
\langle a_{p^r,c}\rangle)}.
\endalign
$$
At this point, we need a sublemma which is completely analogous to
Lemma 4-1. 

\proclaim{Sublemma} Let $G\cdot\langle a\rangle$ be a group generated
by a finite group $G$ and an element $a$ of order a power of $p$ such
that $a$ commutes with any element in $G$ and $G\cap\langle a\rangle=
\langle a^{p^r}\rangle\in G_{(p)}$. Suppose $\langle a\rangle$ acts
trivially on $M$. Then
$$
\chi^{(d)}_p(M;G\cdot\langle a\rangle)=p^{rd}\chi^{(d)}_p(M;G).
\tag5-5
$$
\endproclaim
The proof of this sublemma is analogous to Lemma 4-1. Using this
sublemma and the fact that $(a_{p^r,c})^{p^r}\!\!=c\in C_G(c)$, we see
that the exponent of the previous expression is equal to 
$$
\align
\sum_{[c]_p}\chi^{(d-1)}_p(M^{\langle c\rangle};C_G(c)
\cdot\langle a_{p^r,c}\rangle)
&=p^{r(d-1)}\sum_{[c]_p}\chi^{(d-1)}_p
\bigl(M^{\langle c\rangle};C_G(c)\bigr)\\
&=p^{r(d-1)}\chi^{(d)}_p(M;G).
\endalign
$$
Thus, renaming $r$ as $\ell_d$, the expression $(*)$ finally becomes 
$$
(*)=\Bigl[\!\!\!\!\!\!\!\!\!\!\!\prod_{\ell_1,\dots,\ell_{d-1},\ell_d\ge0}
\!\!\!\!\!\!\!\!\!\!\!
(1-q^{p^{\ell_1}\cdots p^{\ell_{d-1}}p^{\ell_d}})^{p^{\ell_2}p^{2\ell_3}
\cdots p^{(d-2)\ell_{d-1}} p^{(d-1)\ell_d}}\Bigr]
^{(-1)\chi^{(d)}_p(M:G)},
$$
which is the right hand side of formula (5-4). This completes the
proof. 
\qed
\enddemo

Now letting $M$ be a point and using (2-2) with $K=\Bbb Z_p^d$, we get

\proclaim{Corollary 5-3} Let $G_n=G\wr\frak S_n$ for $n\ge0$. For any
$d\ge0$, we have 
$$
\sum_{n\ge0}q^n\bigl|\text{\rm Hom}(\Bbb Z_p^d,G_n)/G_n\bigr|
=\Bigl[\!\!\!\!\!\!\!\prod_{\ell_1,\dots,\ell_d\ge0}
\!\!\!\!\!\!\!\!\!\
(1-q^{p^{\ell_1}\cdots p^{\ell_d}})^{p^{\ell_2}p^{2\ell_3}
\cdots p^{(d-1)\ell_d}}\Bigr]^{(-1)|\text{Hom}(\Bbb Z_p^d,G)/G|}.
\tag5-6
$$
\endproclaim

In particular, letting $G$ to be the trivial group, we get 

\proclaim{Corollary 5-4} For any $d\ge0$, 
$$
\sum_{n\ge0}q^n\bigl|\text{\rm Hom}
(\Bbb Z_p^d,\frak S_n)/\frak S_n\bigr|
=\Bigl[\!\!\!\!\!\!\!\!\!\prod_{\ell_1,\ell_2,\dots,\ell_d\ge0}
\!\!\!\!\!\!\!\!\!\
(1-q^{p^{\ell_1}p^{\ell_2}\cdots p^{\ell_d}})^{p^{\ell_2}p^{2\ell_3}
\cdots p^{(d-1)\ell_d}}\Bigr]^{(-1)}.
\tag5-7
$$
\endproclaim

Observe that $\bigl|\text{\rm Hom} (\Bbb Z_p^d,\frak S_n)/\frak
S_n\bigr|$ is the number of isomorphism classes of $\Bbb Z^d_p$-sets
of order $n$. Any finite $\Bbb Z_p^d$-sets can be decomposed into
transitive $\Bbb Z_p^d$-sets which must have order powers of $p$. For
any $r\ge0$, isomorphism classes of transitive $\Bbb Z_p^d$-sets of
order $p^r$ are in $1:1$ correspondence with index $p^r$ subgroup of
$\Bbb Z_p^d$, by taking isotropy subgroups. Let $j_{p^r}(\Bbb Z_p^d)$
be the number of index $p^r$ subgroup of $\Bbb Z_p^d$. Note that
$j_{\ell}(\Bbb Z_p^d)$ is zero unless $\ell$ is a power of $p$. This
consideration of decomposing finite $\Bbb Z_p^d$-sets into transitive
ones immediately gives the following formula.
$$
\sum_{n\ge0}q^n\bigl|\text{\rm Hom}
(\Bbb Z_p^d,\frak S_n)/\frak S_n\bigr|
=\Bigl[\prod_{r\ge0}(1-q^{p^r})^{j_{p^r}(\Bbb Z_p^d)}\Bigr]^{(-1)}.
\tag5-8
$$
There is an easy way to calculate the number $j_{p^r}(\Bbb Z_p^d)$. 

\proclaim{Lemma 5-5} For any $r\ge0$ and $d\ge1$, we have 
$$
\aligned
j_{p^r}(\Bbb Z_p^d)&=\!\!\!\!\sum_{\sum_i\ell_i=r}\!\!\!\!
p^{\ell_2}p^{2\ell_3}\cdots 
p^{(d-1)\ell_d}=\!\!\!\!\!\!\!\sum_{n_1\cdots n_d=p^r}\!\!\!\!\!\!\!
n_2n_3^2\cdots n_d^{d-1},\\
j_{p^r}(\Bbb Z^d_p)&=\!\!\sum_{0\le\ell\le r}
\!\!p^{\ell}\cdot j_{p^{\ell}}(\Bbb Z^{d-1}_p).
\endaligned
\tag5-9
$$
\endproclaim
\demo{Proof} Any subgroup $H$ of $\Bbb Z_p^d$ of index $p^r$ for any
$r\ge0$ is a closed subgroup and hence it has a structure of a $\Bbb
Z_p$-submodule, and as such it is a free module. Let the standard
basis of $\Bbb Z_p^d$ be $e_1,e_2,\dots,e_d$. It is easy to see that
$H$ has a unique $\Bbb Z_p$-module basis $\{x_i\}_{i=1}^n$ of the form
$x_i=p^{\ell_i}e_i+\sum_{i<j\le d}a_{ij}e_j$ for some uniquely
determined integers $\ell_1,\ell_2,\dots,\ell_d\ge0$ with
$\sum_i\ell_i=r$ and $0\le a_{ij}<p^{\ell_j}$. Any choice of such
integers gives rise to a subgroup of $\Bbb Z_p^d$ of index
$p^r$. Thus, counting all possible choices of these integers subject
to $\sum_i\ell_i=r$, we obtain the expression of $j_{p^r}(\Bbb Z_p^d)$
given in (5-9). 

The second formula in (5-9) is straightforward from the first. This
completes the proof.
\qed
\enddemo

Note that formulae (5-8) and (5-9) give an alternate proof of (5-7). 

Now, formulas (1-16) and (1-17) in the introduction easily follow from
(5-9) by induction on $d\ge1$. In terms of these numbers $j_{p^r}(\Bbb
Z_p^d)$, the formula (5-4) can be rewritten as
$$
\sum_{n\ge0}q^n\chi^{(d)}_p(M^n;G\wr\frak S_n)
=\Bigl[\prod_{r\ge0}(1-q^{p^r})^{j_{p^r}(\Bbb Z_p^d)}\Bigr]
^{(-1)\chi^{(d)}_p(M;G)}.
\tag5-10
$$
Since the Euler characteristic of the equivariant Morava K-theory
$K(d)^*_{G_n}(M^n)$ is equal to the $d$-th order $p$-primary orbifold
Euler characteristic $\chi^{(d)}_p(M^n;G_n)$ by Proposition 5-1, we
obtain our final result of this paper:

\proclaim{Theorem 5-6} For any $d\ge0$ and for any $G$-manifold $M$, 
$$
\sum_{n\ge0}q^n\chi\bigl(K(d)^*_{G_n}(M^n)\bigr)
=\Bigl[\prod_{r\ge0}(1-q^{p^r})^{j_{p^r}(\Bbb Z^d_p)}\Bigr]^{(-1)
\chi(K(d)^*_G(M))},
\tag5-11
$$
where $G_n=G\wr\frak S_n$ is a wreath product. When $M$ is a point,
this formula gives 
$$
\sum_{n\ge0}q^n\chi\bigl(K(d)^*(BG_n)\bigr)
=\Bigl[\prod_{r\ge0}(1-q^{p^r})^{j_{p^r}(\Bbb Z^d_p)}\Bigr]^{(-1)
\chi(K(d)^*BG)},
\tag5-12
$$
where $BG$ and $BG_n$ are classifying spaces of $G$ and $G_n$. 
\endproclaim

\Addresses\recd

\bye

%% file: gtmacros.tex
%
%
%
%
%
%
\magnification=\magstephalf      
%
%
\vsize=7.5truein                 
\hsize=5.2truein                 
\newskip\stdskip                 
\stdskip=6pt plus3pt minus3pt    
\medskipamount=\stdskip          
\parindent=0pt                   
\parskip=\stdskip                
\abovedisplayskip=\stdskip       
\belowdisplayskip=\stdskip       
\mathsurround=0.75pt             
\overfullrule=0pt                
%
%
\def\ppar{\par\goodbreak\vskip 8pt plus 4pt minus 4pt}     
%
%
\def\stdspace{\hskip 0.75em plus 0.15em\ignorespaces}
\let\qua\stdspace 
%
%
%
%
%
%
%
\def\hexnumber#1{\ifcase#1 0\or 1\or 2\or 3\or 4\or 5\or 6\or 7\or 8\or
 9\or A\or B\or C\or D\or E\or F\fi}
%
%
\font\thirtnmsa=msam10 scaled 1315    
\font\tenmsa=msam10          \font\ninemsa=msam9
\font\sevenmsa=msam7         \font\sixmsa=msam6
\font\fivemsa=msam5
%
%
\newfam\msafam                  \textfont\msafam=\tenmsa
\scriptfont\msafam=\sevenmsa    \scriptscriptfont\msafam=\fivemsa
\edef\hexa{\hexnumber\msafam}        
\def\msa{\fam\msafam\tenmsa}         
%
%
\font\thirtnmsb=msbm10 scaled 1315   
\font\tenmsb=msbm10      \font\ninemsb=msbm9
\font\sevenmsb=msbm7     \font\sixmsb=msbm6
\font\fivemsb=msbm5
%
\newfam\msbfam                   \textfont\msbfam=\tenmsb       
\scriptfont\msbfam=\sevenmsb     \scriptscriptfont\msbfam=\fivemsb
\edef\hexb{\hexnumber\msbfam}    
\def\msb{\fam\msbfam\tenmsb}     
%
%
\font\thirtneufm=eufm10 scaled 1315   
\font\teneufm=eufm10                 \font\nineeufm=eufm9
\font\seveneufm=eufm7                \font\sixeufm=eufm6
\font\fiveeufm=eufm5
%
\newfam\eufmfam                    \textfont\eufmfam=\teneufm
\scriptfont\eufmfam=\seveneufm     \scriptscriptfont\eufmfam=\fiveeufm
\edef\hexf{\hexnumber\eufmfam}      
\def\frak{\fam\eufmfam\teneufm}     
%
%
%
\font\thirtnrm=cmr10 scaled 1315    
\font\ninerm=cmr9                   \font\sixrm=cmr6   
%
\font\thirtni=cmmi10 scaled 1315    
\font\ninei=cmmi9                   \font\sixi=cmmi6  
%
\font\thirtnsy=cmsy10 scaled 1315   
\font\ninesy=cmsy9                  \font\sixsy=cmsy6  
%
\font\thirtnbf=cmbx10 scaled 1315   
\font\ninebf=cmbx9                  \font\sixbf=cmbx6  
%
%
\font\thirtnex=cmex10 scaled 1315   
\font\nineex=cmex9                  
%
%
\font\thirtnit=cmti10 scaled 1315  
\font\nineit=cmti9                  
%
\font\thirtnsl=cmsl10 scaled 1315  
\font\ninesl=cmsl9                  
%
\font\thirtntt=cmtt10 scaled 1315  
\font\ninett=cmtt9                  
%
%
%
%
\def\small{%
%
%
\textfont0=\ninerm \scriptfont0=\sixrm \scriptscriptfont0=\fiverm
\def\rm{\fam0\ninerm}
%
%
\textfont1=\ninei \scriptfont1=\sixi \scriptscriptfont1=\fivei
%
%
\textfont2=\ninesy \scriptfont2=\sixsy \scriptscriptfont2=\fivesy
%
%
\textfont3=\nineex \scriptfont3=\nineex \scriptscriptfont3=\nineex
%
%
\textfont\bffam=\ninebf \scriptfont\bffam=\sixbf
\scriptscriptfont\bffam=\fivebf \def\bf{\fam\bffam\ninebf}%
%
%
\textfont\itfam=\nineit \def\it{\fam\itfam\nineit}%
\textfont\slfam=\ninesl \def\sl{\fam\slfam\ninesl}%
\textfont\ttfam=\ninett \def\tt{\fam\ttfam\ninett}%
%
%
%
\textfont\msafam=\ninemsa \scriptfont\msafam=\sixmsa
\scriptscriptfont\msafam=\fivemsa \def\msa{\fam\msafam\ninemsa}%
%
%
\textfont\msbfam=\ninemsb \scriptfont\msbfam=\sixmsb
\scriptscriptfont\msbfam=\fivemsb \def\msb{\fam\msbfam\ninemsb}%
%
%
\textfont\eufmfam=\nineeufm  \scriptfont\eufmfam=\sixeufm
\scriptscriptfont\eufmfam=\fiveeufm \def\frak{\fam\eufmfam\nineeufm}%
%
%
%
\normalbaselineskip=11pt%
\setbox\strutbox=\hbox{\vrule height8pt depth3pt width0pt}%
%
%
\normalbaselines\rm
%
%
\stdskip=4pt plus2pt minus2pt    
\medskipamount=\stdskip          
\parskip=\stdskip                
\abovedisplayskip=\stdskip       
\belowdisplayskip=\stdskip       
\def\ppar{\par\goodbreak\vskip 6pt plus 3pt minus 3pt}%
%
%
\def\section##1{\global\advance\sectionnumber by 1
\vskip-\lastskip\penalty-800\vskip 20pt plus10pt minus5pt 
\egroup{\bf\number\sectionnumber\quad##1}\bgroup\small         
\vskip 6pt plus3pt minus3pt
\nobreak\resultnumber=1}
}    
%
\def\beginsmall{\bgroup\small}
\let\endsmall\egroup
%
%
%
%
\def\large{%
\textfont0=\thirtnrm \scriptfont0=\ninerm \scriptscriptfont0=\sevenrm
\def\rm{\fam0\thirtnrm}%
\textfont1=\thirtni \scriptfont1=\ninei \scriptscriptfont1=\seveni
\textfont2=\thirtnsy \scriptfont2=\ninesy \scriptscriptfont2=\sevensy
\textfont3=\thirtnex \scriptfont3=\thirtnex \scriptscriptfont3=\thirtnex
\textfont\bffam=\thirtnbf \scriptfont\bffam=\ninebf
\scriptscriptfont\bffam=\sevenbf \def\bf{\fam\bffam\thirtnbf}%
\textfont\itfam=\thirtnit \def\it{\fam\itfam\thirtnit}%
\textfont\slfam=\thirtnsl \def\sl{\fam\slfam\thirtnsl}%
\textfont\ttfam=\thirtntt \def\tt{\fam\ttfam\thirtntt}%
\textfont\msafam=\thirtnmsa \scriptfont\msafam=\ninemsa
\scriptscriptfont\msafam=\sevenmsa \def\msa{\fam\msafam\thirtnmsa}%
\textfont\msbfam=\thirtnmsb \scriptfont\msbfam=\ninemsb
\scriptscriptfont\msbfam=\sevenmsb \def\msb{\fam\msbfam\thirtnmsb}%
\textfont\eufmfam=\thirtneufm  \scriptfont\eufmfam=\nineeufm
\scriptscriptfont\eufmfam=\seveneufm \def\frak{\fam\eufmfam\teneufm}%
\normalbaselineskip=16pt%
\setbox\strutbox=\hbox{\vrule height11.5pt depth4.5pt width0pt}%
\normalbaselines\rm}%
\let\Large\large   
%
\def\Bbb#1{{\msb#1}}

%

%
\mathchardef\plussquare="0\hexa01
\mathchardef\nge="3\hexb0B
\mathchardef\maltesecross="0\hexa7A
\mathchardef\del="0\hexf01
%
%
%
%
\font\sc=cmcsc10
%
%
%
%
\def\sqr#1#2{{\vcenter{\vbox{\hrule  height.#2truept
	\hbox{\vrule width.#2truept height#1truept 
	\kern#1truept \vrule width.#2truept}
	\hrule height.#2truept}}}}
\def\sq{\sqr55}    
%
%
%
%
\newcount\sectionnumber            
\newcount\resultnumber             
\sectionnumber=0\resultnumber=1    
%
%
%
\def\section#1{\global\advance\sectionnumber by 1
\xdef\nextkey{\number\sectionnumber}
\vskip-\lastskip\penalty-800\vskip 20pt plus10pt minus5pt 
{\large\bf\number\sectionnumber\quad#1}         
\vskip 8pt plus4pt minus4pt
\nobreak\resultnumber=1}      
%
%
%
%
%
\def\sh#1{\vskip-\lastskip\ppar{\bf #1}\par\nobreak\medskip}         
%
%
%
%

%
\def\proc#1{\xdef\nextkey{\number\sectionnumber.\number\resultnumber}%
\vskip-\lastskip\ppar\bf%
\noindent#1\ \number\sectionnumber.\number\resultnumber
\stdspace\sl\global\advance\resultnumber by 1\ignorespaces}
 
%
%
\def\qed{\hfill$\sq$\par\goodbreak\rm}   
%
%
%
%
%
%
%
%
\def\proclaim#1{\vskip-\lastskip\ppar\bf%
\noindent#1\stdspace\sl\ignorespaces} 
\let\endproclaim\endproc
%
%
%
%
\def\rk#1{\vskip-\lastskip\ppar{\bf #1}\stdspace\ignorespaces}                

%
%
%
%
%
%
\def\label{\xdef\nextkey{\number\sectionnumber.\number\resultnumber}%
\number\sectionnumber.\number\resultnumber
\global\advance\resultnumber by 1}
%
%
%
%
%
%
%
%
%
%
%
%
%
%
%
%
\newcount\refnumber              
\refnumber=1                     
\long\def\reflist#1\endreflist{%
\long\def\thereflist{#1}{\def\refkey##1##2\par{\xdef##1{\number\refnumber}%
\global\advance\refnumber by 1}%
\def\key##1##2\par{\expandafter\xdef%
\csname##1\endcsname{\number\refnumber}%
\global\advance\refnumber by 1}#1\par}}
\long\def\references{%
\penalty-800\vskip-\lastskip\vskip 15pt plus10pt minus5pt 
{\large\bf References}\ppar 
{\leftskip=25pt\frenchspacing    
\small\parskip=3pt plus2pt       
\def\refkey##1##2\par{\noindent  
\llap{[##1]\stdspace}\ignorespaces##2\par}         
\def\key##1##2\par{\noindent  
\llap{[\ref{##1}]\stdspace}\ignorespaces##2\par}  
\def\,{\thinspace}\thereflist\par}}
%
%
%
\newcount\footnotenumber         
\footnotenumber=1                
\def\fnote#1{\xdef\nextkey{\number\footnotenumber}%
{\small\ifnum\footnotenumber>9\parindent=14pt%
\else\parindent=10pt\fi\footnote{$^{\number\footnotenumber}$}%
{\hglue-5pt#1}\global\advance\footnotenumber by 1}}
%
%
%
%
%
%
%
\newcount\figurenumber          
\figurenumber=1                 
\def\caption#1{\xdef\nextkey{\number\figurenumber}%
\cl{\small Figure \number\figurenumber: #1}%
\global\advance\figurenumber by 1}
\def\figurelabel{\xdef\nextkey{\number\figurenumber}%
\cl{\small Figure \number\figurenumber}%
\global\advance\figurenumber by 1}
\long\def\figure#1\endfigure{{\xdef\nextkey{\number\figurenumber}%
\let\captiontext\relax\def\caption##1{\xdef\captiontext{##1}}%
\midinsert\cl{\ignorespaces#1\unskip\unskip\unskip\unskip}\vglue6pt\cl{\small 
Figure \number\figurenumber\ifx\captiontext\relax\else: \captiontext
\fi}\endinsert\global\advance\figurenumber by 1}}
%
%
%
%
%
%
%
\def\nextkey{??}   
%
\def\key#1{\expandafter\xdef\csname #1\endcsname{\nextkey}}
\def\ref#1{\expandafter\ifx\csname #1\endcsname\relax
\immediate\write16{Reference {#1} undefined}??\else
\csname #1\endcsname\fi}
%
%
%
%
%
%
%
\newread\gtinfile
\newwrite\gtreffile
\def\useforwardrefs{
\openin\gtinfile\jobname.ref
\ifeof\gtinfile
\closein\gtinfile
\immediate\write16{No file \jobname.ref}
\else
\closein\gtinfile
\input \jobname.ref
\fi
\immediate\openout\gtreffile \jobname.ref
%
%
\def\key##1{{\def\\{\noexpand}%
\expandafter\xdef\csname ##1\endcsname{\nextkey}%
\immediate\write\gtreffile{\\\expandafter\\\def\\\csname ##1\\\endcsname%
{\nextkey}}}}
%
%
\long\def\reflist##1\endreflist{%
\long\def\thereflist{##1}{\def\refkey####1####2\par{\xdef####1{%
\number\refnumber}{\def\\{\noexpand}\immediate\write\gtreffile
{\\\def\\####1{\number\refnumber}}}\global\advance\refnumber by 1}%
\def\key####1####2\par{\expandafter\xdef%
\csname####1\endcsname{\number\refnumber}%
{\def\\{\noexpand}\immediate\write\gtreffile
{\\\expandafter\\\def\\\csname ####1\\\endcsname{\number\refnumber}}}
\global\advance\refnumber by 1}##1\par}}
\long\def\biblio##1\endbiblio{\reflist##1\endreflist\references}%
%
%
\def\numkey##1{{\def\\{\noexpand}%
\xdef##1{\number\sectionnumber.\number\resultnumber}
\immediate\write\gtreffile{\\\def\\##1%
{\number\sectionnumber.\number\resultnumber}}}}
\def\seckey##1{{\def\\{\noexpand}\xdef##1{\number\sectionnumber}
\immediate\write\gtreffile{\\\def\\##1{\number\sectionnumber}}}}
\def\figkey##1{\xdef##1{\number\figurenumber}%
{\def\\{\noexpand}\immediate\write\gtreffile%
{\\\def\\##1{\number\figurenumber}}}
\number\figurenumber\global\advance\figurenumber by 1}
}   
%
%
%
%
\def\figkey#1{\xdef#1{\number\figurenumber}%
\number\figurenumber\global\advance\figurenumber by 1}
\def\fig#1#2\endfig{%
\midinsert\cl{#2}\vglue6pt\cl{\small Figure #1}\endinsert}
\def\newfig{\number\figurenumber\global\advance\figurenumber by 1}
\def\numkey#1{\xdef#1{\number\sectionnumber.\number\resultnumber}}
\def\seckey#1{\xdef#1{\number\sectionnumber}}
%
%
%
%
%
%
%
%
%
\def\verb{\catcode`\"=\active}       
\def\brev{\catcode`\"=12}            
\brev                                
\verb                                
{\obeyspaces\gdef {\ }}              
{\catcode`\`=\active\gdef`{\relax\lq}}
\def"{%
\begingroup\baselineskip=12pt\def\par{\leavevmode\endgraf}%
\tt\obeylines\obeyspaces\parskip=0pt\parindent=0pt%
\catcode`\$=12\catcode`\&=12\catcode`\^=12\catcode`\#=12%
\catcode`\_=12\catcode`\~=12%
\catcode`\{=12\catcode`\}=12\catcode`\%=12\catcode`\\=12%
\catcode`\`=\active\let"\endgroup}
\brev      
%
%
%
%
%
%
\def\item#1{\par\leavevmode\llap{#1\stdspace}%
\ignorespaces}                             
%
%

%
%
\def\np{\vfil\eject}         
\def\nl{\hfil\break}         
\def\cl{\centerline}         
\def\agt{{\mathsurround=0pt\it$\cal A\mskip-.7mu$lgebraic \&\ 
$\cal G\mskip-2mu$eometric $\cal T\!\!$opology}}  
%
%
%

%
%
%
%
%
\def\title#1{\def\thetitle{#1}}
\def\shorttitle#1{\def\theshorttitle{#1}}
\def\author#1{\edef\previousauthors{\theauthors}
 \ifx\theauthors\relax\def\theauthors{#1}\else
 \def\theauthors{\previousauthors\par#1}\fi}

%
\def\address#1{\edef\previousaddresses{\theaddress}
 \ifx\theaddress\relax\def\theaddress{#1}\else
 \def\theaddress{\previousaddresses\par\vskip 2pt\par#1}\fi}
\def\secondaddress#1{\edef\previousaddresses{\theaddress}
 \ifx\theaddress\relax\def\theaddress{#1}\else
 \def\theaddress{\previousaddresses\par{\rm and}\par#1}\fi}   

\def\email#1{\edef\previousemails{\theemail}
 \ifx\theemail\relax\def\theemail{#1}\else
 \def\theemail{\previousemails\hskip 0.75em\relax#1}\fi}
\def\secondemail#1{\edef\previousemails{\theemail}
 \ifx\theemail\relax\def\theemail{#1}\else
 \def\theemail{\previousemails\hskip 0.75em{\rm and}\hskip 0.75em
 \relax#1}\fi}
\def\url#1{\edef\previousurls{\theurl}
 \ifx\theurl\relax\def\theurl{#1}\else
 \def\theurl{\previousurls\hskip 0.75em\relax#1}\fi}
\def\secondurl#1{\edef\previousurls{\theurl}
 \ifx\theurl\relax\def\theurl{#1}\else
 \def\theurl{\previousurls\hskip 0.75em{\rm and}\hskip 0.75em
 \relax#1}\fi}
\long\def\abstract#1\endabstract{\long\def\theabstract{#1}}
\def\primaryclass#1{\def\theprimaryclass{#1}}
\def\secondaryclass#1{\def\thesecondaryclass{#1}}
\def\keywords#1{\def\thekeywords{#1}}
%
%
\let\\\par\let\thetitle\relax\let\theshorttitle\relax
\let\theauthors\relax\let\theshortauthors\relax
\let\theaddress\relax\let\theshortaddress\relax
\let\theemail\relax\let\theurl\relax
\let\theabstract\relax\let\theprimaryclass\relax
\let\thesecondaryclass\relax\let\thekeywords\relax
%
%
%
%
\long\def\maketitlepage{    

\vglue 0.2truein   

%
{\parskip=0pt\leftskip 0pt plus 1fil\def\\{\par\smallskip}{\large
\bf\thetitle}\par\medskip}   

\vglue 0.15truein 

%
{\parskip=0pt\leftskip 0pt plus 1fil\def\\{\par}{\sc\theauthors}
\par\medskip}%
 
\vglue 0.1truein 

%
{\small\parskip=0pt
{\leftskip 0pt plus 1fil\def\\{\par}{\sl\theaddress}\par}
\ifx\theemail\relax\else  
\vglue 5pt \def\\{\stdspace{\rm and}\stdspace} 
\cl{Email:\stdspace\tt\theemail}\fi
\ifx\theurl\relax\else    
\vglue 5pt \def\\{\stdspace{\rm and}\stdspace} 
\cl{URL:\stdspace\tt\theurl}\fi\par}

\vglue 7pt 

{\bf Abstract}

\vglue 5pt

\theabstract

\vglue 7pt 

{\bf AMS Classification numbers}\quad Primary:\quad \theprimaryclass\par

Secondary:\quad \thesecondaryclass

\vglue 5pt 

{\bf Keywords:}\quad \thekeywords

\np  

}    
%
%
\long\def\makeshorttitle{    


%
{\parskip=0pt\leftskip 0pt plus 1fil\def\\{\par\smallskip}{\large
\bf\thetitle}\par\medskip}   

\vglue 0.05truein 

%
{\parskip=0pt\leftskip 0pt plus 1fil\def\\{\par}{\sc\theauthors}
\par\medskip}%
 
\vglue 0.03truein 

%
{\small\parskip=0pt
{\leftskip 0pt plus 1fil\def\\{\par}{\sl\ifx\theshortaddress\relax
\theaddress\else\theshortaddress\fi}\par}
\ifx\theemail\relax\else  
\vglue 5pt \def\\{\stdspace{\rm and}\stdspace} 
\cl{Email:\stdspace\tt\theemail}\fi
\ifx\theurl\relax\else    
\vglue 5pt \def\\{\stdspace{\rm and}\stdspace} 
\cl{URL:\stdspace\tt\theurl}\fi\par}

\vglue 10pt 


{\small\leftskip 25pt\rightskip 25pt{\bf Abstract}\stdspace\theabstract

{\bf AMS Classification}\stdspace\theprimaryclass
\ifx\thesecondaryclass\relax\else; \thesecondaryclass\fi\par
{\bf Keywords}\stdspace \thekeywords\par}
\vglue 7pt
}    
\let\maketitle\makeshorttitle        
%
%

\def\volumenumber#1{\def\thevolumenumber{#1}}
\def\volumeyear#1{\def\thevolumeyear{#1}}
\def\pagenumbers#1#2{\def\startpage{#1}\def\finishpage{#2}}
\def\published#1{\def\publishdate{#1}}
\def\received#1{\def\receiveddate{#1}}
\def\revised#1{\def\reviseddate{#1}}
\let\reviseddate\relax
\volumenumber{X}
\volumeyear{20XX}
\pagenumbers{1}{XXX}
\published{XX Xxxember 20XX}

\long\def\makeagttitle{   
\agt\hfill      
\hbox to 60truept{\vbox to 0pt{\vglue -14truept{\bf [Logo here]}\vss}\hss}
\break
{\small Volume \thevolumenumber\ (\thevolumeyear)
\startpage--\finishpage\nl
Published: \publishdate}

\vglue .2truein

{\parskip=0pt\leftskip 0pt plus 1fil\def\\{\par\smallskip}{\large
\bf\thetitle}\par\medskip}   
\vglue 0.05truein 

%
{\parskip=0pt\leftskip 0pt plus 1fil\def\\{\par}{\sc\theauthors}
\par\medskip}%
 
\vglue 0.03truein 


{\small\leftskip 25truept\rightskip 25truept{\bf Abstract}\stdspace\theabstract

{\bf AMS Classification}\stdspace\theprimaryclass
\ifx\thesecondaryclass\relax\else; \thesecondaryclass\fi\par
{\bf Keywords}\stdspace \thekeywords\par}\vglue 7truept

}   


\def\Addresses{\bigskip
{\small \parskip 0pt \leftskip 0pt \rightskip 0pt plus 1fil \def\\{\par}
\sl\theaddress\par\medskip \rm Email:\stdspace\tt\theemail\par
\ifx\theurl\relax\else\smallskip \rm URL:\stdspace\tt\theurl\par\fi}}

\def\agtart{
\hoffset 14truemm
\voffset 31truemm
\font\phead=cmsl9 scaled 950
\font\pnum=cmbx10 scaled 913
\font\pfoot=cmsl9 scaled 950
\headline{\vbox to 0pt{\vskip -4.5mm\line{\small\phead\ifnum
\count0=\startpage ISSN numbers are printed here
\hfill {\pnum\folio}\else\ifodd\count0\def\\{ }%
\ifx\theshorttitle\relax\thetitle\else\theshorttitle\fi\hfill{\pnum\folio}
\else\def\\{ and }{\pnum\folio}\hfill\ifx\theshortauthors\relax\theauthors
\else\theshortauthors\fi\fi\fi}\vss}}
\footline{\vbox to 0pt{\vglue 0mm\line{\small\pfoot\ifnum\count0=\startpage
Copyright declaration is printed here\hfill\else
\agt, Volume \thevolumenumber\ (\thevolumeyear)\hfill\fi}\vss}}
\let\maketitle\makeagttitle\let\makeshorttitle\makeagttitle}

%% file: agtout.tex

\def\ifplaintex{\expandafter\ifx\csname documentclass\endcsname\relax}

\def\gtp{{\mathsurround=0pt\it $\cal G\mskip-2mu$eometry \&\ 
$\cal T\!\!$opology $\cal P\!$ublications}}  

\def\recd{{\small Received:\qua\receiveddate\ifx\reviseddate\relax
\else\qquad Revised:\qua\reviseddate\fi\par}} 


\def\volumenumber#1{\def\thevolumenumber{#1}}
\def\volumeyear#1{\def\thevolumeyear{#1}}
\def\papernumber#1{\def\thepapernumber{#1}}
\def\pagenumbers#1#2{\def\startpage{#1}\def\finishpage{#2}}
\def\published#1{\def\publishdate{#1}}

\def\received#1{\def\receiveddate{#1}}
\def\revised#1{\def\reviseddate{#1}}
\def\accepted#1{\def\accepteddate{#1}}
\def\asciititle#1{\def\theasciititle{#1}}
\def\covertitle#1{\def\thecovertitle{#1}}

\long\def\asciiabstract#1{\long\def\theasciiabstract{#1}}
\def\asciikeywords#1{\def\theasciikeywords{#1}}


\let\\\par\let\thevolumenumber\relax
\let\thepapernumber\relax\let\thevolumeyear\relax\let\startpage\relax
\let\finishpage\relax\let\publishdate\relax\let\receiveddate\relax
\let\reviseddate\relax\let\accepteddate\relax\let\theasciititle\relax
\let\thecovertitle\relax\let\theasciiauthors\relax
\let\theasciiabstract\relax\let\theasciikeywords\relax

\let\theasciiemail\relax


\ifplaintex
\font\logobig=cmssbx10 scaled 3836
\font\logomed=cmssbx10 scaled 2557
\else
\font\logobig=cmssbx10 scaled 4200
\font\logomed=cmssbx10 scaled 2800
\fi

\long\def\makeagttitle{   
\count0=\startpage
\agt\hfill      
\hbox to 45truept{\vbox to 0pt{\vglue -13truept{\logomed A\kern -.37em{\logobig 
T}\kern -.38em G}\vss}\hss}
\break
{\small Volume \thevolumenumber\ (\thevolumeyear)
\startpage--\finishpage\nl
Published: \publishdate}

\vglue .25truein

{\parskip=0pt\leftskip 0pt plus
1fil\def\\{\par\smallskip}{\Large\bf\thetitle}\par\medskip} \vglue
0.05truein

%
{\parskip=0pt\leftskip 0pt plus 1fil\def\\{\par}{\sc\theauthors}
\par\medskip}%
 
\vglue 0.03truein 


{\small\leftskip 25truept\rightskip 25truept{\bf Abstract}\stdspace\theabstract

{\bf AMS Classification}\stdspace\theprimaryclass
\ifx\thesecondaryclass\relax\else; \thesecondaryclass\fi\par
{\bf Keywords}\stdspace \thekeywords\par}\vglue 7truept

}   

\ifplaintex
\hoffset 14truemm
\voffset 31truemm
\font\phead=cmsl9 scaled 950
\font\pnum=cmbx10 scaled 913
\font\pfoot=cmsl9 scaled 950
\headline{\vbox to 0pt{\vskip -4.5mm\line{\small\phead\ifnum
\count0=\startpage ISSN 1472-2739 (on-line) 1472-2747 (printed)
\hfill {\pnum\folio}\else\ifodd\count0\def\\{ }%
\ifx\theshorttitle\relax\thetitle\else\theshorttitle\fi\hfill{\pnum\folio}
\else\def\\{ and }{\pnum\folio}\hfill\ifx\theshortauthors\relax\theauthors
\else\theshortauthors\fi\fi\fi}\vss}}
\footline{\vbox to 0pt{\vglue 0mm\line{\small\pfoot\ifnum\count0=\startpage
\copyright\ \gtp\hfill\else
\agt, Volume \thevolumenumber\ (\thevolumeyear)\hfill\fi}\vss}}
\else
\headsep 23pt
\footskip 35pt
\hoffset -4truemm
\voffset 12.5truemm
\font\lhead=cmsl9 scaled 1050
\font\lnum=cmbx10 
\font\lfoot=cmsl9 scaled 1050
\makeatletter
\def\@oddhead{{\small\lhead\ifnum\count0=\startpage ISSN 1472-2739 
(on-line) 1472-2747 (printed)\hfill {\lnum\number\count0}\else\ifodd\count0
\def\\{ }\ifx\theshorttitle\relax \thetitle \else\theshorttitle\fi\hfill
{\lnum\number\count0}\else\def\\{ and }{\lnum\number\count0}
\hfill\ifx\theshortauthors\relax 
\theauthors\else\theshortauthors\fi\fi\fi}}\def\@evenhead{\@oddhead}
\def\@oddfoot{\small\lfoot\ifnum\count0=\startpage\copyright\ \gtp\hfill\else
\agt, Volume \thevolumenumber\ (\thevolumeyear)\hfill\fi}
\def\@evenfoot{\@oddfoot}
\makeatother
\fi
\let\maketitlepage\makeagttitle
\let\makeshorttitle\maketitlepage
\let\maketitle\maketitlepage


\newwrite\gtoutfile
\long\gdef\makeheadfile{  
{\def\\{, }\def\s{ }
\immediate\openout\gtoutfile head.xxx
\immediate\write\gtoutfile{To: math@arxiv.org}
\immediate\write\gtoutfile{Subject: put OR rep NNNNN:ppppp}
\immediate\write\gtoutfile{--text follows this line--}
\immediate\write\gtoutfile{Proxy-for: \ifx\theasciiauthors\relax
\theauthors\else\theasciiauthors\fi\s<\ifx\theasciiemail\relax\theemail\else\theasciiemail\fi>}
\immediate\write\gtoutfile{\noexpand\\}
\immediate\write\gtoutfile{Authors: \ifx\theasciiauthors\relax
\theauthors\else\theasciiauthors\fi}
{\def\\{ }\immediate\write\gtoutfile{Title: \ifx\theasciititle\relax
\thetitle\else\theasciititle\fi}}
\immediate\write\gtoutfile{Subj-class: GT or SG, GR etc}
\immediate\write\gtoutfile{MSC-class: \theprimaryclass\ifx\thesecondaryclass\relax\else, \thesecondaryclass\fi}
\immediate\write\gtoutfile{Journal-ref: Algebraic and Geometric Topology \thevolumenumber\s
(\thevolumeyear) \startpage-\finishpage}
\immediate\write\gtoutfile{Comments: Published by Algebraic and
Geometric Topology at}
\immediate\write\gtoutfile{\s\s\s  http://www.maths.warwick.ac.uk/agt/AGTVol\thevolumenumber/agt-\thevolumenumber-\thepapernumber.abs.html}
\immediate\write\gtoutfile{\noexpand\\}
\immediate\write\gtoutfile{}
\ifx\theasciiabstract\relax
\immediate\write\gtoutfile{\theabstract}\else
\immediate\write\gtoutfile{\theasciiabstract}\fi
\immediate\write\gtoutfile{}
\immediate\write\gtoutfile{\noexpand\\}
\immediate\write\gtoutfile{}
\immediate\closeout\gtoutfile}}  

\def\maketitlepage{\makeagttitle\makeheadfile}
\let\makeshorttitle\maketitlepage
\let\maketitle\maketitlepage


\def\ifplaintex{\expandafter\ifx\csname documentclass\endcsname\relax}

\def\gtp{{\mathsurround=0pt\it $\cal G\mskip-2mu$eometry \&\ 
$\cal T\!\!$opology $\cal P\!$ublications}}  

\def\recd{{\small Received:\qua\receiveddate\ifx\reviseddate\relax
\else\qquad Revised:\qua\reviseddate\fi\par}} 


\def\volumenumber#1{\def\thevolumenumber{#1}}
\def\volumeyear#1{\def\thevolumeyear{#1}}
\def\papernumber#1{\def\thepapernumber{#1}}
\def\pagenumbers#1#2{\def\startpage{#1}\def\finishpage{#2}}
\def\published#1{\def\publishdate{#1}}

\def\received#1{\def\receiveddate{#1}}
\def\revised#1{\def\reviseddate{#1}}
\def\accepted#1{\def\accepteddate{#1}}
\def\asciititle#1{\def\theasciititle{#1}}
\def\covertitle#1{\def\thecovertitle{#1}}

\long\def\asciiabstract#1{\long\def\theasciiabstract{#1}}
\def\asciikeywords#1{\def\theasciikeywords{#1}}


\let\\\par\let\thevolumenumber\relax
\let\thepapernumber\relax\let\thevolumeyear\relax\let\startpage\relax
\let\finishpage\relax\let\publishdate\relax\let\receiveddate\relax
\let\reviseddate\relax\let\accepteddate\relax\let\theasciititle\relax
\let\thecovertitle\relax\let\theasciiauthors\relax
\let\theasciiabstract\relax\let\theasciikeywords\relax

\let\theasciiemail\relax


\ifplaintex
\font\logobig=cmssbx10 scaled 3836
\font\logomed=cmssbx10 scaled 2557
\else
\font\logobig=cmssbx10 scaled 4200
\font\logomed=cmssbx10 scaled 2800
\fi

\long\def\makeagttitle{   
\count0=\startpage
\agt\hfill      
\hbox to 45truept{\vbox to 0pt{\vglue -13truept{\logomed A\kern -.37em{\logobig 
T}\kern -.38em G}\vss}\hss}
\break
{\small Volume \thevolumenumber\ (\thevolumeyear)
\startpage--\finishpage\nl
Published: \publishdate}

\vglue .25truein

{\parskip=0pt\leftskip 0pt plus
1fil\def\\{\par\smallskip}{\Large\bf\thetitle}\par\medskip} \vglue
0.05truein

%
{\parskip=0pt\leftskip 0pt plus 1fil\def\\{\par}{\sc\theauthors}
\par\medskip}%
 
\vglue 0.03truein 


{\small\leftskip 25truept\rightskip 25truept{\bf Abstract}\stdspace\theabstract

{\bf AMS Classification}\stdspace\theprimaryclass
\ifx\thesecondaryclass\relax\else; \thesecondaryclass\fi\par
{\bf Keywords}\stdspace \thekeywords\par}\vglue 7truept

}   

\ifplaintex
\hoffset 14truemm
\voffset 31truemm
\font\phead=cmsl9 scaled 950
\font\pnum=cmbx10 scaled 913
\font\pfoot=cmsl9 scaled 950
\headline{\vbox to 0pt{\vskip -4.5mm\line{\small\phead\ifnum
\count0=\startpage ISSN 1472-2739 (on-line) 1472-2747 (printed)
\hfill {\pnum\folio}\else\ifodd\count0\def\\{ }%
\ifx\theshorttitle\relax\thetitle\else\theshorttitle\fi\hfill{\pnum\folio}
\else\def\\{ and }{\pnum\folio}\hfill\ifx\theshortauthors\relax\theauthors
\else\theshortauthors\fi\fi\fi}\vss}}
\footline{\vbox to 0pt{\vglue 0mm\line{\small\pfoot\ifnum\count0=\startpage
\copyright\ \gtp\hfill\else
\agt, Volume \thevolumenumber\ (\thevolumeyear)\hfill\fi}\vss}}
\else
\headsep 23pt
\footskip 35pt
\hoffset -4truemm
\voffset 12.5truemm
\font\lhead=cmsl9 scaled 1050
\font\lnum=cmbx10 
\font\lfoot=cmsl9 scaled 1050
\makeatletter
\def\@oddhead{{\small\lhead\ifnum\count0=\startpage ISSN 1472-2739 
(on-line) 1472-2747 (printed)\hfill {\lnum\number\count0}\else\ifodd\count0
\def\\{ }\ifx\theshorttitle\relax \thetitle \else\theshorttitle\fi\hfill
{\lnum\number\count0}\else\def\\{ and }{\lnum\number\count0}
\hfill\ifx\theshortauthors\relax 
\theauthors\else\theshortauthors\fi\fi\fi}}\def\@evenhead{\@oddhead}
\def\@oddfoot{\small\lfoot\ifnum\count0=\startpage\copyright\ \gtp\hfill\else
\agt, Volume \thevolumenumber\ (\thevolumeyear)\hfill\fi}
\def\@evenfoot{\@oddfoot}
\makeatother
\fi
\let\maketitlepage\makeagttitle
\let\makeshorttitle\maketitlepage
\let\maketitle\maketitlepage


\newwrite\gtoutfile
\long\gdef\makeheadfile{  
{\def\\{, }\def\s{ }
\immediate\openout\gtoutfile head.xxx
\immediate\write\gtoutfile{To: math@arxiv.org}
\immediate\write\gtoutfile{Subject: put OR rep NNNNN:ppppp}
\immediate\write\gtoutfile{--text follows this line--}
\immediate\write\gtoutfile{Proxy-for: \ifx\theasciiauthors\relax
\theauthors\else\theasciiauthors\fi\s<\ifx\theasciiemail\relax\theemail\else\theasciiemail\fi>}
\immediate\write\gtoutfile{\noexpand\\}
\immediate\write\gtoutfile{Authors: \ifx\theasciiauthors\relax
\theauthors\else\theasciiauthors\fi}
{\def\\{ }\immediate\write\gtoutfile{Title: \ifx\theasciititle\relax
\thetitle\else\theasciititle\fi}}
\immediate\write\gtoutfile{Subj-class: GT or SG, GR etc}
\immediate\write\gtoutfile{MSC-class: \theprimaryclass\ifx\thesecondaryclass\relax\else, \thesecondaryclass\fi}
\immediate\write\gtoutfile{Journal-ref: Algebraic and Geometric Topology \thevolumenumber\s
(\thevolumeyear) \startpage-\finishpage}
\immediate\write\gtoutfile{Comments: Published by Algebraic and
Geometric Topology at}
\immediate\write\gtoutfile{\s\s\s  http://www.maths.warwick.ac.uk/agt/AGTVol\thevolumenumber/agt-\thevolumenumber-\thepapernumber.abs.html}
\immediate\write\gtoutfile{\noexpand\\}
\immediate\write\gtoutfile{}
\ifx\theasciiabstract\relax
\immediate\write\gtoutfile{\theabstract}\else
\immediate\write\gtoutfile{\theasciiabstract}\fi
\immediate\write\gtoutfile{}
\immediate\write\gtoutfile{\noexpand\\}
\immediate\write\gtoutfile{}
\immediate\closeout\gtoutfile}}  

\def\maketitlepage{\makeagttitle\makeheadfile}
\let\makeshorttitle\maketitlepage
\let\maketitle\maketitlepage